\documentstyle{amltd}
\begin{document}
\input amssym.def
\input amssym.tex
\annalsline{151}{2000}
%\received{March 17, 1997}
%\revised{July 21, 1999}
\startingpage{1}
\def\bye{\end{document}}
 \font\tenrm=cmr10

%--------------- Author macros ---------------
\def\underset#1#2{{#2}_{#1}}
\def\eqref#1{(\ref{#1})}
\def\ritem#1{\item[{\rm #1}]}
\catcode`\@=11
\font\twelvemsb=msbm10 scaled 1100
\font\tenmsb=msbm10
%\font\ninemsb=msbm7 scaled 1100%msbm9
\font\ninemsb=msbm10 scaled 800
\newfam\msbfam
\textfont\msbfam=\twelvemsb  \scriptfont\msbfam=\ninemsb
  \scriptscriptfont\msbfam=\ninemsb
\def\msb@{\hexnumber@\msbfam}
\def\Bbb{\relax\ifmmode\let\next\Bbb@\else
 \def\next{\errmessage{Use \string\Bbb\space only in math
mode}}\fi\next}
\def\Bbb@#1{{\Bbb@@{#1}}}
\def\Bbb@@#1{\fam\msbfam#1}
\catcode`\@=12

 \catcode`\@=11
\font\twelveeuf=eufm10 scaled 1100
\font\teneuf=eufm10
\font\nineeuf=eufm7 scaled 1100%eufm9
\newfam\euffam
\textfont\euffam=\twelveeuf  \scriptfont\euffam=\teneuf
  \scriptscriptfont\euffam=\nineeuf
\def\euf@{\hexnumber@\euffam}
\def\frak{\relax\ifmmode\let\next\frak@\else
 \def\next{\errmessage{Use \string\frak\space only in math
mode}}\fi\next}
\def\frak@#1{{\frak@@{#1}}}
\def\frak@@#1{\fam\euffam#1}
\catcode`\@=12

%-------------- Author entries --------------------

\title{Some spherical uniqueness theorems\\ for 
multiple trigonometric series} 
\shorttitle{Spherical Uniqueness Theorems} 

 \acknowledgements{The research of both authors 
was partially supported by NSF grant DMS 9707011 and a grant from the Faculty and Development Program
of the College of Liberal Arts and Sciences, De Paul University.\\
\phantom{From*}1991 {\it Mathematics Subject Classification}.
Primary 42B05, 42B99; Secondary 42B08, 42B15, 42A63.  }
 \twoauthors{J.\  Marshall Ash}{Gang Wang}
\institutions{De Paul University,  Chicago, IL\\
{\eightpoint {\it E-mail addresses\/}: mash@math.depaul.edu}\\
 \hglue.97in {\eightpoint gwang@math.depaul.edu}} 

\phantom{more space}

 \bigbreak \centerline{\bf Abstract}
\bigbreak
We prove that if a multiple trigonometric series is spherically Abel 
sum\-mable everywhere to an everywhere finite function $f(x)$ which is 
bounded below by an integrable function, then the series is the Fourier  
series of $f(x)$ if the coefficients of the multiple trigonometric series 
satisfy a mild growth condition. As a consequence, we show that if a 
multiple trigonometric series is spherically convergent everywhere to an 
everywhere finite integrable function $f(x)$, then the series is the Fourier  
series of $f(x)$. We also show that a singleton is a set of uniqueness. These results are generalizations of a recent theorem of 
J. Bourgain and some results of V. Shapiro.

\section{Introduction and summary of results}\label{s:0}

We start with the question of spherical uniqueness of multiple 
trigonometric series for integrable functions under Abel summability. Greek
letters $\xi, \eta, \cdots$ will denote points of the $d$-dimensional
lattice $\Bbb  Z^d$, Roman letters $x, y,\cdots$ points of the 
$d$-dimensional torus $\Bbb {T}^d = [-\pi, \pi)^d$, $\langle \cdot, \cdot
\rangle$ inner product, and $|\cdot|$ $d$-dimensional Euclidean norm.
For a multiple trigonometric series 
$\sum_{\xi\in Z^d} a_{\xi} e^{i\langle x, \xi \rangle}$ 
where the coefficients $a_{\xi}$ are arbitrary complex numbers, the Abel 
sum is defined to be  the limit of the function 
\[
f(x, t) ={\sum}_{\xi\in Z^d} a_{\xi} 
e^{i\langle x, \xi \rangle-|\xi|t}
\]
as $t\to 0^+$ if such limit exists. In general, denote
\begin{eqnarray*}
\noalign{\vskip4pt}
f^*(x) & =& \limsup_{t\to 0^+} f(x,t)\\
\noalign{\vskip4pt}
& =& \Re f^*(x) + i \Im f^*(x) 
\\
\noalign{\vskip-6pt}
\end{eqnarray*}
$f_*(x) = \Re f_*(x) + i \Im f_*(x)$ is similarly defined with $\limsup$ being 
replaced by $\liminf$.

It is well-known, when $d=1$, that if $\sum a_{\xi}e^{i\xi x}$ is Abel 
summable to 0 everywhere and if $a_{\xi}=o(|\xi|)$, then all $a_{\xi}=0$.
See, for example, \cite{Ve1} and \cite{Ve2}. To see that this theorem is
sharp, look at the one dimensional series\break $\delta ^{\prime }(x)=-
\sum \xi\sin {\xi x}$, which may be thought of as the derivative
of the Dirac delta function. It is easy to check that this series is 
Abel summable to $0$, although the growth condition is just 
barely violated. Thinking of $\delta ^{\prime }$ as a 
degenerate $d$ dimensional function, it is immediately clear that the 
hypothesis of a $d$ dimensional uniqueness theorem concerning Abel 
summability will necessarily have to carry some growth condition. 
One generalization of this fact is due to Victor Shapiro, who extended one 
dimensional work of Verblunsky and of Rajchman and Zygmund (\cite{Sh}).

\phantom{hee}

\proclaimtitle{Shapiro}
\proclaim{Theorem} \label{t:sh1}
Let $\sum a_{\xi}e^{i\langle\xi, x\rangle}$ be a multiple trigonometric 
series{\rm .}  Suppose that
\begin{itemize}
\ritem{1.}  the coefficients $a_{\xi}$ satisfy the following growth rate 
condition\/{\rm :}
\end{itemize}
\begin{equation}\label{0.1}
\sum_{R-1<\left| \xi\right| \leq R}\left| a_{\xi}\right| 
=o\left(R\right)  \hbox{ as } R\rightarrow \infty ,
\end{equation}
\begin{itemize}
\ritem{2.}  $f^{*}(x)$ and $f_{*}(x)$ are finite for all $x${\rm , }

\ritem{3.}  $a_{\xi}=\overline{a_{-\xi}}$ for all $\xi${\rm ,} and

\ritem{4.}  $\min\{\Re f_{*}(x), \Im f_{*}(x)\} \geq A(x)$ where $A\left( x\right) 
$ is in $L^{1}\left( \Bbb  {T}^{d}\right) .$
\end{itemize}

\noindent Then $f_{*}\left( x\right) \in L^{1}\left( \Bbb {T}^{d}\right) $ and $\sum
a_{\xi}e^{i\langle \xi,x\rangle}$ is the Fourier series of $f_{*}.$
\endproclaim 

\phantom{hee}

This theorem is sharp because the example $\delta ^{\prime }$
mentioned above just barely fails to meet condition \eqref{0.1}. Nevertheless
condition \eqref{0.1} is disappointingly strong in the sense that when Abel 
summability is replaced by regular convergence, condition \eqref{0.1} is not a direct consequence of convergence. However, there is a known theorem concerning the coefficients' growth 
rate for spherically convergent multiple trigonometric series. In fact, it is 
implied by the following Cantor-Lebesgue type theorem.

\proclaimtitle{Connes}
\proclaim{Theorem} \label{t:con} Let ${\cal O}\subset \Bbb  T^d$
be a ball or a subset which has full measure and is of Baire second category 
relative to $\Bbb  {T}^d$.
If $\sum_{|\xi|=R} a_{\xi} e^{i\langle x, \xi \rangle}$ tends 
to $0$ as $R\to\infty$ at every point of ${\cal O}$, then 
\begin{equation}\label{0.2}
\varepsilon_R^2 = \sum_{|\xi| = R}|a_{\xi}|^2
= o(1) \qquad \hbox{ as } R\to\infty.
\end{equation}
\endproclaim

Connes proved this theorem for dimension $d$ in 1976, twenty years 
after Shapiro's Theorem~\ref{t:sh1}. Cooke \cite{C} and shortly thereafter 
Zygmund \cite{Z} had completed the $d=2$ case five years before Connes' work.

An easy corollary of Theorem~\ref{t:con} gives the coefficients' growth rate
condition for spherically convergent multiple trigonometric series.
 
\proclaimtitle{Connes}
\proclaim{{C}orollary} \label{c:con}
 Let ${\cal O}\subset \Bbb  T^d$
be a ball or a subset which has full measure and is of Baire second category 
relative to $\Bbb  {T}${\rm .}
If\/ $\lim_{R\to\infty}\sum_{|\xi|\leq R} a_{\xi} 
e^{i\langle x, \xi \rangle}$ exists {\rm (}\/as a finite number\/{\rm )} at each point of
${\cal O}${\rm ,} then 
\begin{equation}\label{0.3}
\varepsilon_R^2 = \sum_{|\xi| = R}|a_{\xi}|^2
= o(1) \qquad \hbox{ as } R\to\infty.
\end{equation}
\endproclaim

The coefficients' growth rate condition given by \eqref{0.3} does 
not imply condition \eqref{0.1} when $d\geq 3$. To remedy this problem, we 
first prove the following analogue of Theorem~\ref{t:sh1} under the condition \eqref{0.3}. We use notation $A\sim B$ to denote $B/2 \leq A < B$. 

\proclaim{Theorem} \label{t:1} Consider the multiple trigonometric series 
$\sum_{\xi\in Z^d} a_{\xi} e^{i\langle x, \xi \rangle}$ where
the coefficients $a_{\xi}$ are arbitrary complex numbers{\rm .} Suppose that
\begin{itemize}
\ritem{1.}  the coefficients of the series $a_{\xi}$ satisfy
\end{itemize}
\begin{equation}\label{e:bou}
\sum_{|\xi|\sim R} |a_{\xi}|^2 = \underset{R/2 \leq |\xi|< R}
{\sum} |a_{\xi}|^2 = o(R^2) \hbox{ as } R\to\infty,
\end{equation}
\begin{itemize}
\ritem{2.}   $f^*(x)$ and $f_*(x)$ are finite for all $x${\rm ,} and
\ritem{3.}  $\min\{\Re f_*(x), \Im f_*(x)\}$ is bounded below by a function $A(x)$ 
in $L^1(\Bbb  T^d)${\rm .}
\end{itemize}
\noindent Then $f_*(x)$ is in $L^1(\Bbb  T^d)$ and $\sum_{\xi\in Z^d} 
a_{\xi} e^{i\langle x, \xi \rangle}$ is its Fourier series{\rm .}
\endproclaim

Note that condition \eqref{0.3} implies condition \eqref{e:bou} since
\[
\sum_{|\xi|\sim R} |a_{\xi}|^2 = \sum_{k= R^2/4}^{R^2 - 1} \sum_{|\xi|^2 = k} 
|a_{\xi}|^2 = o\left(\sum_{k= R^2/4}^{R^2 - 1} 1 \right) = o(R^2).
\]

Since \eqref{0.1} implies 
\[
\sum_{|\xi|\sim R} |a_{\xi}| = o(R^2),
\]
\eqref{0.1} implies \eqref{e:bou} if $a_{\xi}$ is bounded. But in 
general, and when $d\geq 3$, \eqref{0.1} and \eqref{e:bou} do not relate to 
each other. Notice that \eqref{e:bou} implies that 
\[
\sum_{|\xi|^2\sim R}|a_{\xi}|^2 = \sum_{R/2 \leq |\xi|^2<R}|a_{\xi}|^2 \leq \sum_{|\xi|\sim \sqrt{R}}|a_{\xi}|^2 = o(R).
\]

As a consequence of Theorem~\ref{t:con} and Theorem~\ref{t:1}, we obtain the following two
spherical uniqueness theorems for multiple trigonometric series which are convergent to a function. These theorems make no assumption whatsoever about coefficient size.

\proclaim{Theorem}\label{t:2}
Let $\sum_{\xi\in Z^d} a_{\xi} e^{i\langle x, \xi \rangle}$ be 
a trigonometric series which converges spherically everywhere to an 
everywhere finite function $f(x)${\rm ;} i{\rm .}e{\rm .,}
\begin{equation}\label{0.4}
\lim_{R\to\infty}\sum_{|\xi|\leq R} a_{\xi} e^{i\langle x, \xi\rangle} = f(x) 
\hbox{ for all } x\in \Bbb {T}^d.
\end{equation}
If $\min\{\Re f(x), \Im f(x)\}\geq g(x)$ for all $x$ and 
$g(x) \in L^1(\Bbb {T}^d)${\rm ,} then
$f(x)$ is in $L^1(\Bbb {T}^d)$ and  $a_{\xi} $ is the $\xi^{\rm th}$ Fourier 
coefficient  of $f(x)$ for all $\xi\in \Bbb {Z}^d${\rm . } 
\endproclaim

In particular,

\proclaim{Theorem}\label{t:3}
Let $f(x)\in L^1(\Bbb {T}^d)$ be finite at every $x${\rm .} If 
$\sum_{\xi\in Z^d} a_{\xi} e^{i\langle x, \xi \rangle}$ is 
a trigonometric series which converges spherically to $f(x)$ at every 
point $x${\rm ,} i{\rm .}e{\rm . }
\begin{equation}\label{0.5}
\lim_{R\to\infty}\sum_{|\xi|\leq R} a_{\xi} e^{i\langle x, \xi\rangle} = f(x) 
\hbox{ for all } x\in \Bbb {T}^d,
\end{equation}
then $a_{\xi} $ is the $\xi^{\rm th}$ Fourier coefficient  of $f(x)$ for all 
$\xi\in \Bbb {Z}^d${\rm .  }
\endproclaim

Special cases of Theorem~\ref{t:3} have been proved by various people. 
When $d=1$ and $f(x) \equiv 0$, this is the original uniqueness theorem of 
Cantor. For general $f(x)\in L^1(\Bbb  T^1)$, it was first proved by de 
la Vall\'ee-Poussin. When  $d=2$, Theorem~\ref{t:sh1} combined with
the work of Cooke \cite{C} implies Theorem~\ref{t:1} and thus, 
Theorem~\ref{t:2} and Theorem~\ref{t:3}.  The major breakthrough came
when Bourgain \cite{B} proved Theorem~\ref{t:3} for the special case of
$f(x) \equiv 0$. 
For a survey on the uniqueness of multiple trigonometric series under 
various summation modes,  as well as many open problems in this area, please 
refer to Ash and Wang \cite{AW}. 

The proof of Theorem~\ref{t:1} is mainly based on Shapiro's framework 
\cite{Sh}. To avoid assuming condition \eqref{0.1}, we exploit an idea that 
Bourgain \cite{B} 
used when he proved Theorem~\ref{t:3} for the special case $f(x) \equiv 0$. 
We refer to \eqref{e:bou} hereafter as Bourgain's condition, in his honor. This  condition simply asserts that Connes'
condition holds ``on the average." 

The detailed proof of Theorem~\ref{t:1} is given in Sections~\ref{s:1}
through \ref{s:4}. 

At the end of the paper, we begin the study of sets of uniqueness for spherical convergence. As a first step toward establishing this theory, we show that any singleton is a set of uniqueness.  
 
\proclaim{Theorem}\label{c:0.1}
Let $q$ be a point on $\Bbb  T^d${\rm .} Suppose that a multiple trigonometric 
series $\sum_{\xi\in Z^d} a_{\xi} e^{i\langle x, \xi \rangle}$
spherically converges everywhere except at $q$ to a function $f(x)\in L^1 (\Bbb  {T}^d)${\rm .}
Furthermore{\rm ,} suppose $f(x)$ is finite for all $x$ except $q${\rm .} Then 
$\sum_{\xi\in Z^d} a_{\xi} e^{i\langle x, \xi \rangle}$ is 
the Fourier series of $f(x)${\rm .}
\endproclaim

It is easily deduced from Theorem~\ref{t:con} and the following fact about Abel summability,  
which is an analogue of a theorem of Shapiro \cite[\S 6]{Sh}:   
  
\proclaim{Theorem}\label{t:4}
Consider the multiple ($d\geq 2$) trigonometric series\break 
$\sum_{\xi\in Z^d} a_{\xi} 
e^{i\langle x, \xi \rangle}$ where the coefficients $a_{\xi}$ are arbitrary
complex numbers{\rm .}  Let $q$ be a point on $\Bbb  T^d${\rm .} Suppose that
\begin{itemize}
\ritem{1.}  $\sum_{|\xi|\sim R} |a_{\xi}|^2 = o(R^2)$ as
$R\to\infty${\rm ,}
\ritem{2.}  $f^*(x)$ and $f_*(x)$ are finite for all $x$ except $q${\rm ,} and
\ritem{3.}  $f^*(x)$ and $f_*(x)$ are functions in $L^1(\Bbb  T^d)${\rm .}
\end{itemize}
\noindent Then $\sum_{\xi\in Z^d} a_{\xi} e^{i\langle x, \xi \rangle}$ is 
the Fourier series of $f_*(x)${\rm .}
\endproclaim

Note that the theorem is false when $d=1$ since the trigonometric series $\sum
e^{i\xi x}$ is Abel convergent to $0$ everywhere in $\Bbb  T \setminus
\{0\}$. 
 	
\section{Proof of Theorem~\ref{t:1}}\label{s:1}

We may assume $d\ge 3$ since the cases $d=1$ and $d=2$ are known. 

We need some preliminary results and some notation before we start the 
proof. Without loss of generality, by considering the real and imaginary 
parts  separately, we may assume that $a_{\xi} = \overline a_{-\xi}$, where 
$\overline a$ is the conjugate of the complex number $a$. Thus $f(x, t), 
f^*(x)$ and $f_*(x)$ are all real functions. In addition, we may assume 
that  $a_0 = 0$. 

Define 
\begin{eqnarray}\label{1.0}
f_1(x,t) &=& -\underset{\xi \neq 0}{\sum} \frac{a_{\xi}}{|\xi|} 
e^{i\langle x, \xi \rangle-|\xi|t}\\
\noalign{\vskip6pt}
f_2(x,t) &=& -\underset{\xi \neq 0}{\sum} \frac{a_{\xi}}{|\xi|^2} 
e^{i\langle x, \xi \rangle-|\xi|t}.
\nonumber
\end{eqnarray}  
Under the condition \eqref{e:bou}, it is easy to see that for each 
$x\in \Bbb {T}^d$ and $t>0$, $f(x, t), f_1(x, t)$ and $f_2(x, t)$ converge 
absolutely and hence are infinitely differentiable as functions of $t>0$. 
Thus by the mean value theorem, for $t_1>t_2>0$, there exist $t_3, t_4\in 
(t_2, t_1)$ such that $f_1(x, t_1)-f_1(x, t_2) = f(x, t_3)(t_1-t_2)$, and 
$f_2(x, t_1)- f_2(x, t_2) = -f_1(x, t_4)(t_1-t_2)$. Since for each $x$, 
$f^*(x)$ and $f_*(x)$ are finite,  $f(x, t)$ is bounded for all $t>0$. The 
bound depends on $x$ in general. Thus, for each $x$, there exist finite-valued functions $f_1(x)$ and $f_2(x)$ such that   
\begin{equation}\label{1.001}
f_1(x, t)\to f_1(x) \hbox{ and } f_2(x, t)\to f_2(x) \hbox{ as } t\to 0^+.
\end{equation}

On the other hand, if we define the Riemann function $F(x)$ by
\begin{equation}\label{1.1}
F(x) = -\sum_{\xi \neq 0} \frac{a_{\xi}}{|\xi|^2} e^{i\langle x, \xi\rangle},
\end{equation}
then, because of the Bourgain condition \eqref{e:bou}, $F(x)\in L^2(\Bbb  T^d)$ 
and $f_2(x, t)\break \stackrel{L^2}{\rightarrow} F(x)$ as $t\to 0^{+}$. In fact,  
observe that there is an absolute constant $C$ such that 
\begin{eqnarray}\label{1.01}
\qquad ||f_2(x, t)-F(x)||_2^2 &=& \sum_{\xi\neq
0}\frac{|a_{\xi}|^2}{|\xi|^4}\left(1-e^{-2|\xi|t}\right)\\
& \leq & C \sum_{k = 1}^{\infty} 2^{-2k}\left(1 - e^{-2^{k/2+1}t}\right)
\sum_{|\xi|^2 \sim 2^{k}}|a_{\xi}|^2\nonumber \\
&\leq& C \sum_{k=1}^{\infty} 2^{-k}\left(1- e ^{-2^{k/2+1}t}\right)\nonumber\\
&\to& 0 \hbox{ as } t\to 0^+. \nonumber
\end{eqnarray} 
Thus, 
\begin{equation}\label{1.2}
f_2(x) =F(x) \quad {\rm a.e.}
\end{equation}

The key to the proof is to show that $\Delta f_2(x) = f_*(x)$ almost 
everywhere. To this end, we need to use a generalized Laplacian. 

Let $B(x, \rho)$ be an open ball in $\Bbb {T}^d$ centered at $x\in 
\Bbb {T}^d$
with radius $\rho >0$ and $m(B(x, \rho))$ the volume of $B(x,
\rho)$. Then  $m(B(x, \rho)) = v_d \rho^d$, where $v_d$ is the volume of
the unit ball in $\Bbb  R^d$. For any locally integrable function $g(x)$, the average 
of $g$ over $B(x, \rho)$ is
\begin{eqnarray*}
A_\rho g(x) &= &\frac{1}{m(B(x, \rho))}\int_{B(x, \rho)} g(y) \, dy\\
\noalign{\vskip4pt}
            & =& \frac{1}{v_d\rho^d}\int_{B(x, \rho)} g(y) \, dy.
\end{eqnarray*}

Let
\[
I(x) = \frac{I_{B(0, 1)}(x)}{m(B(0, 1))},
\]
where $I_{B(0, 1)}(x)$ is the characteristic function of the unit ball. 
Denote
$\hat I(\xi)$ to be the Fourier transform of $I(x)$. Then $\hat I(\rho
\xi)$ satisfies the following properties:
\begin{equation}\label{1.3}
\lim_{\rho\to 0} \frac{\hat I(\rho\xi)-1}{\rho^2 |\xi|^2} = -\frac12
\int_{B(0, 1)} x_1^2 I(x) \, dx = c_d < 0,
\end{equation}
and for $|\xi| = 1$,
\begin{equation}\label{1.4}
\int_0^\infty \left|\partial_r\left[\frac{\hat I(\rho r \xi) -1}{\rho^2
r^2}\right]\right| \, dr = \int_0^\infty \left|\partial_r\left[\frac{\hat
I(r \xi) -1}{r^2}\right]\right| \, dr < c.
\end{equation}
Note that the constant $c$ in \eqref{1.4} is independent of $\rho$.

The above two equalities are standard. In fact, to see \eqref{1.3}, rotate
 (choose the first coordinate axis to be in the direction of 
$\xi$) and use polar coordinates to get 
\begin{eqnarray}\label{1.5}
\hat I(\rho \xi) -1 &= &\frac{1}{v_d} \int_{|x|\le 1}
(e^{i\langle x, \rho\xi\rangle}-1)\, dx   \\
\noalign{\vskip4pt}
& =& \frac{1}{v_d} \int_{|x| \le 1} (e^{i
x_1\rho|\xi|}-1)\, dx\nonumber\\ \noalign{\vskip4pt}
& = &\frac{v_{d-1}}{v_d} \int_{-1}^{1} (\cos(\rho |\xi| x_1)-1)
(1-x_1^2)^{\frac{d-1}{2}}\, dx_1.\nonumber
\end{eqnarray}
Since for any $x\in \Bbb {T}$, $|\cos x -1| \le \frac{x^2}{2}$, and 
$\lim_{x\to 0} \frac{\cos x -1}{x^2} = -\frac12$,
by the bounded convergence theorem and \eqref{1.5}, we have
\begin{eqnarray*}
\lim_{\rho\to 0} \frac{\hat I(\rho\xi)-1}{\rho^2 |\xi|^2} &=&
-\frac{v_{d-1}}{2v_d} \int_{-1}^{1} x_1^2 (1-x_1^2)^{\frac{d-1}{2}}\, dx_1\\ \noalign{\vskip4pt}
& =& -\frac{1}{2}\int_{B(0,1)} x_1^2 I(x) \, dx = c_d < 0.
\end{eqnarray*}
Observe that the above argument shows $\hat I(\xi_1) = \hat
I(\xi_2)$ if $|\xi_1| = |\xi_2|$. Thus, we may abuse our notation and write 
$\hat I(\xi) = \hat I(|\xi|)$.

Inequality \eqref{1.4} also follows similarly. If $|\xi| = 1$, then
\[
\partial_r \left[\frac{\hat I(\rho r \xi) -1}{\rho^2 r^2}\right] =
-2\frac{v_{d-1}}{v_d} 
\int_{-1}^{1} \frac{(1- \frac{1}{2}i r \rho x_1) e^{i r \rho x_1} -1}{r^3
\rho^2} (1-x_1^2)^{\frac{d-1}{2}}\, dx_1.
\]
Thus, for $c= -2v_{d-1}/v_{d}$, 
\begin{eqnarray*}
&& 
\hskip-36pt\int_0^\infty \left|\partial_r\left[\frac{\hat I(\rho r \xi) -1}{\rho^2
r^2}\right]\right| \, dr\\
  &&\qquad \enspace = \  c \int_{0}^{\infty}\left| \int_{-1}^{1}
\frac{(1- \frac{1}{2}i r \rho x_1) e^{i r \rho x_1} -1}{r^3 \rho^2}
(1-x_1^2)^{\frac{d-1}{2}}\, dx_1 \right|\, dr\\  \noalign{\vskip4pt}  
&&\qquad \enspace = \   c \int_{0}^{\infty}\left| \int_{-1}^{1}
\frac{(1 -\frac{1}{2} i r x_1) e^{i r x_1} -1}{r^3 }
(1-x_1^2)^{\frac{d-1}{2}}\, dx_1 \right|\, dr\\ \noalign{\vskip4pt} 
&&\qquad \enspace = \   \int_0^\infty \left|\partial_r\left[\frac{\hat
I(r \xi) -1}{r^2}\right]\right| \, dr
\end{eqnarray*}
by the simple change of variable argument: $\rho r \to r$.
The above integral is finite since
\[
\left|\int_{-1}^{1} \frac{(1- \frac{1}{2}i r x_1) e^{i r x_1} -1}{r^3
}(1-x_1^2)^{\frac{d-1}{2}} \, dx_1 \right| \le
\frac{c}{r^2} 
\hbox{ as } r \to \infty ,
\]
and
\[
(1- \frac{1}{2} i r x_1) e^{i r x_1} -1 = \frac{1}{2} i r x_1 + O(r^3)
\hbox{ as } r \to 0.
\]

Define the generalized Laplacian  operator  on $g(x)\in L^1$ to be
\[
\tilde \Delta g(x) = \lim_{\rho\to 0}-\frac{1}{c_d} 
\frac{A_\rho g(x)-g(x)}{\rho^2}
\]
if such a limit exists (not necessarily finite), where $c_d<0$ is the
constant given in \eqref{1.3}. We can also define the upper and lower 
generalized Laplacians $\tilde\Delta^{*} g(x)$ and $\tilde\Delta_{*} g(x)$ 
by replacing $\lim$ by $\limsup$ and $\liminf$ respectively  when the 
function $g(x)$ is real-valued. It is clear that all three of these generalized Laplacians agree with the usual Laplacian when applied to a $C^2$ function. Recall that $a_0 = 0$. For $f_2(x, t)$ given 
by \eqref{1.0}, we have for $x\in \Bbb {T}^d$,  
\begin{eqnarray}\label{1.6}
\frac{A_{\rho} f_2(x, t) -f_2(x, t)}{\rho^2} & = & - \sum_{\xi\ne 0} 
\frac{a_{\xi}}{|\xi|^2}
\frac{\hat I(\rho \xi) -1}{\rho^2} e^{i\langle x, \xi\rangle-|\xi|t}\\  \noalign{\vskip4pt}
& = &- \sum_{k\ge 1} \frac{\hat I(\rho \sqrt{k}) -1}{\rho^2 k}\sum_{|\xi|^2 =
k} a_{\xi} e^{i\langle x, \xi \rangle-|\xi|t}\nonumber\\ \noalign{\vskip4pt} 
& =&- \sum_{k\ge 1} \left(\sum_{|\xi|^2 \le k} a_{\xi}e^{i\langle x,\xi
\rangle-|\xi|t}\right)\nonumber\\ \noalign{\vskip4pt} 
&& \qquad\quad\times \left(\frac{\hat I(\rho\sqrt{k}) -1}{\rho^2 k}-\frac{\hat
I(\rho\sqrt{k+1}) -1}{\rho^2 (k+1)}\right)\nonumber \\ \noalign{\vskip4pt} 
& \to & - f(x, t) c_d  \hbox{ as } \rho \to 0 \nonumber
\end{eqnarray}
since by the fundamental theorem of calculus and \eqref{1.4},
\[
\sum_{k\ge 1} \left|\frac{\hat I(\rho\sqrt{k}) -1}{\rho^2 k} -\frac{\hat
I(\rho\sqrt{k+1}) -1}{\rho^2 (k+1)}\right| \le \int_0^\infty
\left|\partial_r\left[\frac{\hat I(\rho r \xi) -1}{\rho^2
r^2}\right]\right| \, dr < b <\infty 
\]
for a constant $b$ independent of $\rho$.  Thus, the above argument shows 
that  
\begin{equation}\label{1.7}
\tilde \Delta f_2(x, t) = \lim_{\rho\to 0}-\frac{1}{c_d}
\frac{A_{\rho} f_2(x, t)-f_2(x, t)}{\rho^2}  =
f(x, t) 
\end{equation}
for $x\in \Bbb {T}^d$ and $t> 0$.

To pass to the limit as $t\to 0^+$, we need the following lemma of Shapiro 
(Lemma~7 of \cite{Sh2}).  To see that Shapiro's lemma applies,
note that $F \in L^2(\Bbb {T}^d)$ implies $F \in L^1(\Bbb {T}^d)$.
\proclaim{Lemma}\label{l:1.0}
If 
{\rm \eqref{e:bou}} holds{\rm ,} then at every point $x$ where $f_*(x)$ and $f^*(x)$ are 
finite{\rm ,}  
\begin{equation}\label{1.8}
\tilde\Delta_* f_2(x) \leq f^*(x) \quad \hbox{ and } \quad f_*(x) 
\leq \tilde\Delta^* f_2(x).
\end{equation}
\endproclaim

The following classical results on the Green's function $G(x)$ appear with proof as  
Lemma 8 of Shapiro \cite{Sh2}. (Also see Theorem 6 of Bochner \cite{Bo}.)

\proclaim{Lemma}\label{l:1.1}
There is a function $G(x)$ in $L^1(\Bbb {T}^d)$
whose Fourier series is given by $\sum_{|\xi|\ne 0} |\xi|^{-2}
e^{i\langle x, \xi\rangle}${\rm .} Further{\rm ,} $G(x)$ has the following properties\/{\rm :}
\begin{itemize}
\ritem{1.} $G(x)$ is in class $C^{\infty}(\Bbb {T}^d)$ away from $0$ and $\Delta G(x) = 1$ for $x\ne 0${\rm .} 
\ritem{2.} $G(x) = \Phi (x) + H^{*}(x)$ where $H^{*}$ is continuous on
$\Bbb  T^d$ and $\Delta H^*(x) = 1$ for $x\in\Bbb {T}^d\setminus \{0\}$ and where $\Phi(x) = C_d |x|^{-(d-2)}$ 
for $d\ge 3$ with 
$C_d = 2^{d-1}{\pi}^{d/2}\Gamma(d/2)/(d-2)$ and $\Phi(x)=-2\pi\log |x|$ when $d=2$. 
\ritem{3.}  Let $u(x)$ be an upper semi\/{\rm -}\/continuous function on $\Bbb  T^d$
which is also\break in $L^1(\Bbb  T^d)${\rm .} Define $U(x) = (2\pi)^{-d}\int_{T^d} 
G(x-y) u(y)\, dy$ and $u_0 =\break (2\pi)^{-d}\int_{T^d} u(y)\, dy${\rm .}  Then $U(x)$
is upper semi\/{\rm -}\/continuous on $\Bbb {T}^d${\rm ,}\break $U(x)\in L^1(\Bbb  T^d)${\rm ,} and 
$\tilde\Delta_{*}U(x) \ge -u(x) + u_0$ for $x\in \Bbb {T}^d${\rm .} Moreover{\rm ,}
$\tilde\Delta^{*}U(x) = \tilde\Delta_{*}U(x) = - u(x) + u_0$ almost 
everywhere in $\Bbb  T^d${\rm .} 
\end{itemize}

\endproclaim

A consequence of Lemma~\ref{l:1.1} is that for any integrable  
function $u$, the Fourier series of $U(x) = (2\pi)^{-d}\int_{T^d} G(x-y)
u(y)\, dy$ is $\sum_{|\xi|\ne 0} u_{\xi} |\xi|^{-2}e^{i\langle x, \xi
\rangle}$, where $u_0 + \sum_{|\xi|\ne 0}  u_{\xi}e^{i\langle x, \xi
\rangle}$ is the Fourier series of $u$. 

We now state the following key lemma which will be proved in 
Section~\ref{s:2}. The function $\overline U$ will not in general be periodic, so we have to work in $\Bbb  R^d$, rather than in $\Bbb  T^d$. 

\proclaim{Lemma}\label{l:1.2}
Let $f_2(x)$ be as  given in {\rm \eqref{1.001}} where $f(x, t)$ satisfies the 
conditions in Theorem~{\rm \ref{t:1}.} Suppose that $\overline U(x)$ is an 
upper semi\/{\rm -}\/continuous function and that it is in $L_{\rm loc}^1(\Bbb  {R}^d)${\rm .}
 Let $S(x) = f_2(x)+\overline U(x)${\rm .} If 
$\tilde\Delta^* S(x) \geq 0${\rm ,} then $S(x)$ is subharmonic in $\Bbb  R^d${\rm .}
\endproclaim

\demo{{R}emark\/ {\rm 2.1}}
By modifying the proof of Lemma~\ref{l:1.2} in Section~\ref{s:2}, 
Lemma~\ref{l:1.2} can be shown to hold locally. Explicitly, we can replace 
$\Bbb 
R^d$ everywhere in Lemma~\ref{l:1.2} by any open ball $B\subset \Bbb 
R^d$ and $L_{\rm loc}^1(\Bbb  {R}^d)$ by $L^1(B)$.
\enddemo

We now are ready to prove Theorem~\ref{t:1}. 

Since $A(x)\in L^1(\Bbb {T}^d)$, there exists an upper 
semi-continuous function $u(x)$ (see p.75 of \cite{S}, for example) such 
that $u(x) \leq A(x)$. As in Lemma~\ref{l:1.1}, define 
$U(x) = (2\pi)^{-d}\int_{T^d} G(x-y) u(y)\, dy$,  
$u_0 = (2\pi)^{-d}\int_{T^d} u(y)\, dy$ and 
$S(x) = f_2(x) + U(x) - u_0|x|^2/(2d)$. Then by Lemma~\ref{l:1.0}, 
$\tilde\Delta^{*} f_2(x) \geq f_*(x) \geq A(x)  \geq u(x)$. Consequently, 
by periodicity, Lemmas~\ref{l:1.1} and \ref{l:1.2}, $S(x)$ is subharmonic in 
$\Bbb  R^d$. Therefore, by Riesz's representation for subharmonic functions and a theorem of Saks \cite{S}, $\tilde\Delta^{*} S(x) = \tilde\Delta_{*} S(x)$ almost everywhere and is in $L^{1}$ locally. Since $\tilde\Delta^* U(x) = \tilde\Delta_* U(x)$ almost everywhere and is in $L^{1}$ locally, this shows that $\tilde\Delta^{*} f_2(x) = \tilde\Delta_{*} f_2(x)$ almost everywhere and is in $L^{1}$ locally. Thus by assumption and Lemma~\ref{l:1.0}, $f_*(x)$ is in $L^1$ locally.  

Let $B(x) = \min\{f^*(x), \tilde\Delta^{*} f_2(x)\}$. Then by Lemma~\ref{l:1.0}, $\tilde\Delta_* f_2(x) \leq B(x) \leq \tilde\Delta^{*} 
f_2(x)$. Consequently,  $B(x) = \tilde\Delta^{*} f_2(x)$ almost everywhere, 
and is in $L_{\rm loc}^{1}(\Bbb  R^d)$. By a theorem of Vitali-Carath\'eodory (p. 75 of \cite{S}), there exists a nondecreasing sequence of upper semi-continuous functions $\{u^k(x)\}$ on $\Bbb {R}^d$, which are also in $L_{\rm loc}^{1}(\Bbb  R^d)$, such that each $u^k(x)$
 is bounded above and $u^{k}(x) \leq B(x)$ for all $x\in \Bbb  R^d$, 
\begin{equation}\label{1.9}
\lim_{k\to\infty} u^k(x) = B(x) \hbox{ for almost all } x\in \Bbb {R}^d
\end{equation}
and 
\begin{equation}\label{1.10}
\lim_{k\to\infty}\int_E u^{k}(y)\, dy= \int_E B(y)\,dy
\end{equation}
for any bounded set $E\subset\Bbb {R}^d$. 
Set $U^k(x) = (2\pi)^{-d}\int_{T^d} G(x-y) u^k(y)\, dy$ and   
$u^k_0 = (2\pi)^{-d}\int_{T^d} u^k(y)\, dy$. Then \eqref{1.10} implies that 
$u^k_0$ is convergent to\break $b_0 = (2\pi)^{-d}\int_{T^d} B(y)\, dy$ as 
$k \to \infty$. By Lemma~\ref{l:1.0}--Lemma~\ref{l:1.2},
we have $S^k(x) =  f_2(x) + U^k(x) -u^k_0|x|^2/(2d)$ is subharmonic in
$\Bbb {R}^d$.  

Note that $0\leq  B(x) - u^k(x) \leq  B(x) - u^1(x)$. Since $B(x)$ and $u^1(x)$ 
are locally integrable  on $\Bbb {R}^d$, by Lemma~\ref{l:1.1}, \eqref{1.9}, and
the dominated convergence theorem, 
\[
\lim_{k\to\infty} U^{k}(x) =  U(x) = (2\pi)^{-d}\int_{T^d} G(x-y) B(y)
\,  dy  \quad\hbox{ in } L_{\rm loc}^1(\Bbb  R^d)
\]
and hence there exists a subsequence, still called $U^k$ for notational
simplicity, such that 
\[
\lim_{k\to\infty} U^k(x) = U(x)\quad\hbox{ a.e.}
\]
Since for any sequence of subharmonic functions convergent in $L^1$, there is a subharmonic function which is almost everywhere the $L^1$ limit of that sequence (see p. 20 of \cite{R}); $S(x) = f_2(x) + U(x)-b_0|x|^2/(2d)$ is almost everywhere equal to a subharmonic function $S_*(x)$ in $\Bbb  R^d$.

Similarly, there exists a sequence of nonincreasing lower semi-continuous
functions $v^k(x)$ on $\Bbb {R}^d$, which are also
in $L_{\rm loc}^1(\Bbb {R}^d)$, such that each $v^{k}(x)$ is bounded below and $v^{k}(x) \geq B(x)$, 
\[
\lim_{k\to\infty} v^k(x) = B(x) \hbox{ for almost all } x\in \Bbb {R}^d
\]
and 
\[
\lim_{k\to\infty}\int_E v^{k}(y)\, dy= \int_E B(y)\,dy
\]
for any bounded set $E\subset \Bbb {R}^d$. Since $-v^k(x)$ is 
nondecreasing the above arguments show that there exists a superharmonic 
function $S^{*}(x)$, which is almost  everywhere  equal to $S(x)$.

Therefore $S_*(x) = S^*(x)$ almost  everywhere. The subharmonicity of $S_*$
and superharmonicity of $S^*$ show that at every $x$
\begin{equation}\label{1.11}
S_*(x) \leq  A_1 S_*(x) = A_1 S(x) = A_1 S^*(x) \leq  S^*(x).
\end{equation}
In addition, if both $S_*(x)$ and $S^*(x)$ are finite, for any $\epsilon >0$, there exists $\delta>0$, such that
\begin{equation}\label{1.12}
S_*(y) \leq  S_*(x) + \epsilon \hbox{ and } S^*(y) \geq S^{*}(x) -\epsilon
\end{equation}
for all $y\in B(x, \delta)$. Thus the fact that $S_* = S^*$ almost  everywhere and \eqref{1.12} imply that $ S^*(x) \leq  S_*(x) + 2\epsilon$.
So $S^*(x) \leq  S_{*}(x)$. In fact, a similar argument shows that for all $x$, $S_*(x) > -\infty$ and $S^*(x) < \infty$ since $S_*(x) < \infty$ and $S^*(x) > - \infty$ by sub- or superharmonicity. Thus $S_*(x)$ and $S^*(x)$ are finite for all $x$ and $S^*(x) \leq S_*(x)$. Consequently, by \eqref{1.11} $S^*(x) = S_*(x)$ everywhere and hence it is harmonic in $\Bbb {R}^d$.

But then,
\begin{eqnarray*}
\noalign{\vskip4pt}
S^*(x) & =& \frac{1}{v_d}\int_{|y-x|\leq  1} S^*(y)\, dy\\ \noalign{\vskip4pt}
       & =& \frac{1}{v_d}\int_{|y-x|\leq  1} S(y)\, dy\\\noalign{\vskip4pt}
     & =& \frac{1}{v_d}\int_{|y-x|\leq  1} f_2(y)\, dy +
\frac{1}{v_d}\int_{|y-x|\leq  1} U(y)\, dy \\\noalign{\vskip4pt}
&& \qquad-\frac{b_0}{2d}
\frac{1}{v_d}\int_{|y-x|\leq  1} |y|^2\, dy\\\noalign{\vskip4pt}
     & = &I + II + III.
\end{eqnarray*}
Since $f_2$ and $U$ are periodic, $I$ and $II$ are bounded. Thus, $S^*(x) =
O(|x|^2)$. By the penultimate inequality of Section 2.13 of \cite{PW} it follows that every second order partial derivative of $S^{*}$ is a bounded harmonic function and hence constant, so that $S^{*}$ itself is a quadratic polynomial. (An alternative argument can be based on expanding $S^*(x)$ into spherical harmonics.) Thus, the periodic function $f_{2} +U$ is almost everywhere equal to a quadratic polynomial $Q(x)=c_{1,0,\cdots,0}x^2_{1}+\cdots$. A simple countability argument shows that for almost every $x\in \Bbb  R^d$ we have $(f_2+U)(x+2\pi ne_{1})=Q(x+2\pi ne_{1})$ for $n=1,2,3,\cdots$, where $e_1= (1,0,\cdots,0).$ Let $n\to\infty$ to see that $c_{1,0,\cdots,0}=0$. Similar reasoning shows that $Q(x)$ reduces to a constant $K$. 
Consequently, we have $f_2(x) = -U(x) +K$ almost  everywhere. However,
both $U$ and $F$ are integrable over $\Bbb  T^d$. The integrals of $U$ and
$f_2$ over $\Bbb  T^d$ are both $0$ by \eqref{1.2} and Lemma~\ref{l:1.1}.
So $K = 0$. Hence, 
\begin{eqnarray*}
\noalign{\vskip-9pt}
f_2(x) & = &  - (2\pi)^{-d}\int_{T^d}G(x-y)B(y)\, dy\quad {\rm a.e.} \\
& =& -(2\pi)^{-d}\int_{T^d}G(x-y)\tilde\Delta^*f_2(y)\, dy\quad {\rm a.e.}\\
\noalign{\vskip-16pt}
\end{eqnarray*}
Finally,  by \eqref{1.2}, we have 
\[
F(x) = -(2\pi)^{-d}\int_{T^d}G(x-y)\tilde\Delta^* f_2(y)\, dy\quad {\rm a.e.}
\]
Comparing the Fourier series of both sides, we see that the $a_{\xi}$ 
are the Fourier 
coefficients of $\tilde\Delta^* f_2(x)- K_1$ for some constant $K_1$. The
Fourier series of the integrable function $\tilde\Delta^* f_2(x)- K_1$, $\sum
a_{\xi} e^{i\langle \xi, x\rangle}$ is Abel summable to $\tilde\Delta^*
f_2(x)- K_1$ almost everywhere (Theorem 2 of \cite{Sh2}). Thus, from the
definition of $f_*(x)$,  $f_*(x) = \tilde\Delta^* f_2(x)- K_1$ almost
everywhere. Therefore $f_*(x) \in L^1(\Bbb {T}^d)$ and $a_{\xi}$ is the 
$\xi^{\rm th}$ Fourier coefficient of $f_*(x)$. In fact, $K_1 = 0$ by 
Lemma~\ref{l:1.0}.  This  completes the proof that 
Lemma~\ref{l:1.2} will imply Theorem~\ref{t:1}. 
 
We end this section with the following observation. It is well known that if $u(x)$ is an 
upper semi-continuous function in $B= B(x_0, h_0)\subset \Bbb  T^d$ and $G_B$ denotes the Green function of $B$, 
then when $d\geq 3$, the function 
\[
U'(x) = \frac{1}{\sigma_d(d-2)}\int_{B} G_B(x, y) u(y)\, dy ,
\]
where $\sigma_d$ is the surface area of the unit ball in $\Bbb  R^d$,
satisfies
\[
\tilde\Delta_{*}  U'(x) \geq - u(x), \hbox{ for all } x\in B(x_0, h_0).
\]
Replacing $U$ everywhere by $U'$ in the proof of Theorem~\ref{t:1}, we
have the following lemma. Notice that we include the case where $a_0$ may not be zero.  

\proclaim{Lemma}\label{l:1.3} Let 
$\sum_{\xi\in Z^d} a_{\xi} e^{i\langle x, \xi \rangle}$ be a 
multiple $(d\geq 3)$ trigonometric series with
 $\overline a_{\xi} = a_{-\xi}${\rm .} Suppose that the coefficients  $a_{\xi}$ 
satisfy condition  {\rm \eqref{e:bou}; }
\begin{itemize} 
\ritem{1.} $f^*(x)$ and $f_*(x)$ are finite for all $x\in B$ where $B\subset
\Bbb  T^d$ is a ball\/{\rm ;} and
\ritem{2.} $f_*(x) \geq A(x)$ for almost all $x\in B${\rm ,}  where
$A(x)$ is in  $L^1(B)${\rm .}
\end{itemize}
\noindent Then for any ball $B_1\subset \overline{B_1}\subsetneqq B${\rm ,} $f_*$ is in $L^1(B_1)${\rm .} Moreover{\rm ,  }
\[
\overline f_2(x) = f_2(x)+
\frac{1}{\sigma_d(d-2)}\int_{B_1}G_{B_1}(x, y) f_*(y)\, dy
+ a_0|x|^2/(2d)
\]
is finite everywhere and is almost everywhere equal to a function $h(x)$
harmonic in $B_1${\rm .} In addition{\rm ,} if $f^*(x) = f_*(x)$
everywhere in $B${\rm ,} is in $L^1(B)${\rm ,}  and  is continuous in $B${\rm ,} then 
\[
f_2(x) + \frac{1}{\sigma_d(d-2)}\int_{B}G_{B}(x, y) f_*(y)\, dy
+ a_0|x|^2/(2d)
\]
is harmonic on $B${\rm .}
\endproclaim

Note that under condition \eqref{e:bou}, 
\[
F(x) = -\sum_{\xi\neq 0}\frac{a_{\xi}}{|\xi|^2}e^{i\langle \xi, x\rangle}
\]
is in $L^{2}(\Bbb  T^d)$ and $F(x) = f_2(x)$ almost everywhere in $\Bbb 
T^d$. Combining the above lemma with Lemma~5 of Shapiro \cite{Sh1}, we have
the following analogue of Lemma~3 of Shapiro \cite{Sh4}.

Let $B^{o}$ denote the interior of $B$.

\proclaim{Lemma}\label{l:1.4} Let 
$\sum_{\xi\in Z^d} a_{\xi} e^{i\langle x, \xi \rangle}$
be a multiple trigonometric series 
with
$\overline a_{\xi} = a_{-\xi}${\rm .} Suppose that the coefficients $a_{\xi}$
satisfy condition {\rm \eqref{e:bou},}
\begin{itemize} 
\ritem{1.} $f^*(x)$ and $f_*(x)$ are finite for all $x\in B$ where $B\subset
\Bbb  T^d$ is a ball {\rm (}\/open or closed\/{\rm ),} and
\ritem{2.} $f_*(x) =0$ for almost all $x\in B${\rm .}
\end{itemize}
\noindent Then for any ball $B_1 \subset \overline{B_1}\subset B^{o}$, $f(x,t)$ converges to $0$ as 
$t
\to 0^+$ uniformly in $B_1${\rm .} In particular{\rm ,} $f^*(x) = f_*(x) = 0$ in
$B${\rm . }
\endproclaim

\section{Proof of Lemma~\ref{l:1.2}}\label{s:2}

The proof of Lemma~\ref{l:1.2} is so difficult that this section will be given the following preface.

The proof of Lemma~\ref{l:1.2} is extremely delicate, incorporating all the subtle ideas from Bourgain's landmark work \cite{B} as well as an additional Baire category argument that overcomes the unpleasant fact that an upper semi-continuous function on a compact set
need  not be uniformly upper semi-continuous. Some of the difficulty is pushed into Lemmas~\ref{l:2.1} and
\ref{l:2.2}. The proof of Lemma~\ref{l:2.1} contains a great deal of hard analysis. Even after   so much of the
work in Lemmas~\ref{l:2.1} and \ref{l:2.2} has been hidden, the reasoning involved in the proof of Lemma~\ref{l:1.2} is
still tortuous, and so we will provide an overview here.

We assume that the set $W$ where $S$ fails to be upper semi-continuous is nonempty and then reason down a path which eventually divides into two paths each ending in a contradiction. First, a Baire category argument produces a nonempty portion $Z$ of $W$ ($Z=W\cap B$ for some ball $B$) such that $S$ restricted to $Z$ is ``very good," $f_2$ restricted to $Z$ is ``very good,"
et~cetera. 

Next, for each $\varepsilon>0$, let $W_\varepsilon$ be the points of $W$ where $S$ has a jump of at least~$\varepsilon$: 
\[
\limsup_{y\to x} S(y)\geq S(x)+\varepsilon, \hbox{ for every } x\in W_\varepsilon.
\]
For each $\varepsilon>0$ and each $x\in B\setminus \overline W_{\varepsilon}$ consider the harmonic measure $\omega$ of $\partial W_\varepsilon$ with respect to $B\setminus \overline W_\varepsilon$ at $x$. Our path splits depending on whether the harmonic measure is ``thin:"
\begin{equation}\label{I}
\omega (B\setminus\overline W_\varepsilon, \partial W_\varepsilon, x) = 0
\hbox{ for all pairs } (\varepsilon, x) \hbox{ with } \varepsilon >0 \hbox{ and } x\in B\setminus\overline W_\varepsilon,
\end{equation}
or whether it is ``thick:"
\begin{equation}\label{II}
\omega (B\setminus\overline W_\varepsilon, \partial W_\varepsilon, x) > 0
\hbox{ for some } \varepsilon >0 \hbox{ and some } x\in B\setminus\overline W_\varepsilon.
\end{equation}

If \eqref{I} holds, from Lemma~\ref{l:2.1} it follows that $S$ is bounded above and Lemma~\ref{l:2.2} then applies and asserts that $W\cap B=\emptyset$, a contradiction.

On the other hand, if \eqref{II} is the case, we apply a second Baire category argument to strengthen assumption \eqref{II} by producing an $\varepsilon>0$ and a subset $Z_\varepsilon$ of $W_\varepsilon$ so that $U$ is ``uniformly $\varepsilon/40$ subharmonic" when restricted to $Z_\varepsilon$. Furthermore, the set $Z_\varepsilon$ is still ``thick:" $\omega (B\setminus\overline Z_\varepsilon, \partial Z_\varepsilon, x) >0$. 

Finally a very careful procedure involving picking balls within balls within balls 
is used to find a point $p_1$ of $W_\varepsilon$ and a very nearby point $p_2$ of $B$ so that $S(p_2)-S(p_1)$ is small
relative to $\varepsilon$ because of Lemma~\ref{l:2.1}, but large relative to $\varepsilon$ because of $S$ having large
(relative to $\varepsilon$) jumps at each point of $W_\varepsilon$. This contradiction will complete the proof of
Lemma~\ref{l:1.2} which we begin here.

Since $\tilde\Delta^* S(x) \geq 0$ and $S(x)$ is in $L^1$ locally, we have  $S(x)<\infty$ for all $x\in \Bbb  R^d$.
$S(x)\not\equiv -\infty$ since $S(x)$ is in $L^1$ locally. 
We first show that $S(x) = f_2(x) + \overline U(x)$ is 
upper semi-continuous in $\Bbb {R}^d$. 

Let
\begin{equation}\label{2.1}
W_{\varepsilon} = \left\{x\in \Bbb  R^d: \sup_{|x-y| <\delta}
 S(y)-S(x) > \varepsilon \hbox{
for all } \delta > 0\right\}.
\end{equation}
Then the set where $S(x)$ in $\Bbb {R}^d$ is not upper semi-continuous is
given by
\[
W =\bigcup_{\varepsilon >0} W_\varepsilon .
\]

If $W = \emptyset$, then $S(x)$ is upper-semicontinuous. Now we assume 
$W\ne \emptyset$ and  construct the set $Z$. Bourgain's condition \eqref{e:bou} implies that $f(x,t)=\sum a_{\xi} e^{i\langle \xi,x\rangle-|\xi|t}$ is a uniform limit of its partial sums and hence is continuous on $\Bbb {T}^d\times [\frac1j,\infty)$ for every positive integer $j$. Taking periodicity into account, we see that for each $k$, $f(x,t)$ is uniformly continuous on $\Bbb {R}^d\times [\frac{1}{k+1},\frac1k]$. So we may partition $[\frac12,1]$ into $1=t_1>t_2>\cdots>t_r=\frac12$ so that for $i=1,2,\cdots,r-1$, 
\begin{equation}\label{e:locdelta}
\sup_{ x\in R^d}\,\sup_{t_i\ge t\ge t_{i+1}} |f(x,t_i)-f(x,t)|\le 1\ .
\end{equation}
Then partition $[\frac13,\frac12]$ into $\frac12=t_r>t_{r+1}>\cdots>t_s=\frac13$ so that inequality \eqref{e:locdelta} holds for $i=r,r+1,\cdots,s-1$ and so on, thereby producing a sequence 
${\cal T} = \{t_{n}\}$ satisfying $1=t_1>t_2>\cdots$, $\lim_{k\to\infty} t_k =0$, and 	
\begin{equation}\label{e:locdelta1}
\sup_{ x\in R^d}\,\sup_{t_i\ge t\ge t_{i+1}} |f(x,t_i)-f(x,t)|\le 1 
\end{equation}
holds for all $k$. Since $f(x,t)$ is bounded as $t\to 0^+$ for each $x\in\Bbb  R^d$,  
\[
\bigcup_{n\geq 1}\bigcap_{t\in {\cal T}}\left\{x\in\Bbb 
R^d: |f(x, t)|\leq n \right\} = 
\Bbb  R^d.
\]
Therefore,
\[
\bigcup_{n\geq 1}\bigcap_{t\in {\cal T}} \left\{x\in 
\overline W: |f(x, t)|\leq n \right\} 
= \overline W. 
\]
Since for each positive integer $n$ and each $t\in{\cal T}$, the set
\[
\left\{x\in \overline W:
|f(x, t)|\leq n \right\}
\]
is relatively closed with respect to $\overline W$, by Baire's
category theorem applied to the space $\overline W$ (the intersection of countably many relatively open dense sets
 is not empty), for some $N_0\geq 1$, 
\[
\bigcap_{t\in {\cal T}}\left\{x\in \overline W: |f(x, t)| \leq  N_0  \right\}
\]
has a nonempty interior relative to $\overline W$. This means that there exist an
open ball $B(p, \rho_0), p\in W$, and a constant $N_0$ such that
\begin{equation}\label{2.06}
\sup_{t\in {\cal T}} \,\sup_{x\in B(p, \rho_0)\cap\overline W} \,
|f(x, t)|\leq N_0 <\infty.
\end{equation}
Bourgain's condition implies that $\sup_{x\in\Bbb  R^d} |f(x,t)|\le C$, $\sup_{x\in\Bbb  R^d}|f_1(x,t)|\le C$, and $\sup_{x\in\Bbb  R^d}|f_2(x,t)|\le C$ whenever $t\ge 1$. Use 
\[
f_1(x,t)=\int_{1}^{t} f(x,s)\, ds +f_1(x,1),
\] 
\eqref{e:locdelta1} and \eqref{2.06} to see that there is a constant $N>0$ such that 
\begin{equation}\label{l:2.01-1}
\sup_{x\in Z\atop t>0}  |f_1 (x,t)|\le N,
\end{equation}
where $Z=B(p, \rho_0)\cap\overline W$. 
Similarly, since  
\[
f_2(x,t)=-\int_{0}^{t} f_1(x,s)\, ds +f_2(x)
\] 
by \eqref{1.001},  
\begin{equation}\label{2.3}
\sup_{{x\in Z\atop t>0}}|f_2(x) - f_2(x, t)|\leq N t.
\end{equation}
Therefore, $f_2(x)$ is continuous when restricted to $Z$. It follows that $S(x)$ is upper semi-continuous restricted to $Z$.

We will show  a contradiction if $W\ne \emptyset$. Once $S$
is everywhere upper semi-continuous, it is subharmonic since $\tilde\Delta_* S\geq 0$. For this, see p.14 of \cite{R}. This
will complete the proof of Lemma~\ref{l:1.2}.

The following lemmas are needed in proving $W = \emptyset$.  

\smallbreak

For a bounded open set $G$ and Borel measurable set $F$, we denote $\omega (G, F, x)$ to be the harmonic measure of a Borel
set
$F$ relative to $G$ at $x\in G$. Harmonic measure is closely related to Brownian motion. Let $(\{X_t\}_t,  {\cal F}_t, P)$ 
be the standard Brownian motion in $\Bbb {R}^d$.  For $x\in G$, let $T$ be the exiting time of $X_t$ from~$G$:
\[
T = \inf\{t\geq 0: X_t\notin G\}.
\] 
Then $X_T\in \partial G$ since $X_t$ is continuous in $t$. Let $P^x$
denote the probability measure such that $X_0 = x$ almost  everywhere. 
Then the harmonic measure $\omega (G, F, x) = P^x(X_T\in F)$. 

The following properties of harmonic measure are well-known. We summarize 
them as a preliminary lemma.

\proclaim{Lemma}\label{l:2.0} Let $F_0\subset F_1\subset F_2$ be closed subsets of a
bounded open set~$G${\rm .} Then for $x\in G\setminus F_2${\rm ,}
\begin{eqnarray}
\quad \omega (G\setminus F_2, \partial F_2, x)&\geq &\omega (G\setminus F_1, \partial F_1, x)\\
\noalign{\vskip4pt}
&\geq&
\omega (G\setminus F_1, \partial F_0, x)\geq\omega (G\setminus F_2, \partial F_0, x).
\nonumber
\end{eqnarray}
\endproclaim 
 
To see the last inequality, let $T_i = \inf\{t\geq 0: X_t\notin G\setminus F_i\},\, i = 1, 2$. Then $T_2\leq  T_1$. Note
that  on
$\{ T_2< T_1\}$, we must have $X_{T_2} \in \overline G\setminus F_1$. Otherwise, $X_{T_2}\in \overline G\setminus
(\overline G \setminus F_1) = F_1$. Thus by definition $T_2 \geq T_1$, a contradiction. But $\{T_2< T_1\}\subset \{
X_{T_2}\in
\overline G\setminus F_1\}$ implies that $\{ X_{T_2}\in F_1\} \subset \{
T_1 = T_2\}$. Consequently, $\{ X_{T_2}\in \partial F_0\}\subset \{
X_{T_2}\in F_0\}\subset \{
X_{T_2}\in F_1\} \subset \{ T_1 = T_2\}$. This proves the inequality since
\[
P^x(X_{T_2}\in \partial F_0) = P^x(X_{T_2}\in \partial F_0, T_1 = T_2) \leq 
P^x(X_{T_1}\in \partial F_0).
\]

The middle inequality is simply the monotonicity of the harmonic measure. 
To see the left inequality, observe that  on $\{ X_{T_1}\in \partial F_1\}$, $X_{T_2} \in \partial F_2$. Otherwise,
$X_{T_2}\in
\partial G$ and $T_2\leq T_1$ imply that $X_{T_1}\in \partial G$, a contradiction. Consequently, $\{ X_{T_1}\in \partial
F_1\} \subset \{ X_{T_2}\in \partial F_2\}$. This completes the proof.

The next three lemmas are essential to the proof of Lemma~\ref{l:1.2}.

\proclaim{Lemma}\label{l:2.1} Let $S${\rm ,} $f_2${\rm ,} and $\overline U$ be as given
in Lemma~{\rm \ref{l:1.2}.} Let $W$ be the set where $S$ is not upper
semi\/{\rm -}\/continuous{\rm .}
 Assume that there is a open ball $B(p, \rho_0)${\rm ,} $p\in W${\rm ,}
 such that when restricted to $Z =B(p, \rho_0)\cap\overline W${\rm ,}
$f_2(x)$ is continuous and {\rm \eqref{2.3}} holds{\rm .} Then{\rm ,}
 for $p_1\in W, B(p_1,\rho_1) \subset B(p, \frac12\rho_0)$ and $p_2\in
B(p_1, \frac12\rho_1)${\rm ,} there exists a constant $c>0$ such that for almost all such $\rho_1${\rm ,} 
\begin{eqnarray}\label{2.4}
\qquad \enspace S(p_2) - S(p_1)  &\leq & c\biggl( \left[|f_2(p_1)|
+\rho_1^{-\frac34(d-1)}\right]\\ \noalign{\vskip4pt}
 && \times\ [1-\omega (B(p_1, \rho_1)\setminus\overline
W, \partial(W\cap B(p_1, \rho_1)), p_2)]^{\frac14}\nonumber \\ \noalign{\vskip4pt}
  && +\ \sup_{q\in B(p_1, 2\rho_1)\cap\overline W} \, |f_2(q)- f_2(p_1)|\nonumber\\ \noalign{\vskip4pt}
&& +\ 2\sup_{q\in B(p_1, 2\rho_1)}\, \left(\overline U(q) 
-\overline U(p_1)\right)\biggr). \nonumber
\end{eqnarray}
\endproclaim

Lemma~\ref{l:2.1} is a one-sided version of Bourgain's key lemma in 
\cite{B}.  The proof is also similar and is given in Section~\ref{s:3}. It
follows from Lemma~\ref{l:2.1} that $S(x)$ is bounded from above in $B(p,
\frac{\rho_0}{4})$ when  $p_1 =p$.

\proclaim{Lemma}\label{l:2.2} Assume  $\overline U$ is  defined on 
$\overline B(p, r)$ and is upper semi\/{\rm -}\/continuous on  $B(p, r)${\rm .} Let $f_2$ be 
a function in $\overline B(p, r)$ such 
that $S(x)= f_2(x) + \overline U(x)$ is bounded from above in $\overline B(p,
r)${\rm ,} in $L^1(\overline B(p, r))${\rm ,}  and satisfies
\begin{equation}\label{2.5}
\tilde\Delta^* S(x)  \geq 0, \hbox{ and } \tilde\Delta_* f_2(x) <\infty
\end{equation}
for each $x\in B(p, r).$
If $S$ is upper semi\/{\rm -}\/continuous when restricted to 
$\overline W= \{x \in B(p, r): $S(x)$ \hbox{ is not upper 
semi\/{\rm -}\/continuous}\}${\rm ,} and for all $x\in B(p, r)\setminus\overline W_{\varepsilon}$ the harmonic measure
\begin{equation}\label{2.7}
\omega(B(p, r)\setminus\overline W_{\varepsilon}, B(p,r)\cap W_{\varepsilon}, x) = 0 
\end{equation}
for all $\varepsilon>0$ where $W_{\varepsilon}$ is given by {\rm \eqref{2.1},}
 then $W$ must be empty and $S(x)$ is subharmonic on $B(p, r)${\rm . }
\endproclaim

The proof of Lemma~\ref{l:2.2} is given in Section~\ref{s:4}. The special
case when  $\overline U\equiv 0$ was proved by Bourgain. 

The next lemma provides a harmonic measure version of a point density.

\proclaim{Lemma}\label{l:2.3}  Let $B(p_0, r)$ be a ball in $\Bbb  R^d$ and
$F$ a closed set such that $B(p_0, r)\cap F\ne\emptyset${\rm .} Suppose for some 
$x\in B(p_0, r)\setminus F${\rm , }
\[
\omega(B(p_0, r)\setminus F, \partial (B(p_0, r)\cap F), x) >0.
\]
Then there exists $p_1\in B(p_0, r)\cap F${\rm ,} such that
\begin{equation}\label{2.8}
\inf_{\delta_1> 0}\liminf_{\delta_2\to 0}\inf_{x\in B(p_1,
\delta_2)} \, \omega(B(p_1, \delta_1)\setminus F, \partial (B(p_1,
\delta_1)\cap F), x) = 1.
\end{equation}
\endproclaim	

The proof of Lemma~\ref{l:2.3} is outlined in \cite{B}. For a detailed proof,
see the proof of Theorem 3.14 in \cite{AW}. 

We now return to the proof of Lemma~\ref{l:1.2}. 
By \eqref{1.2} and the fact that 
$F(x) \in L^2(\Bbb  T^d)$,  we have $S$ is in $L^1_{\rm loc}(\Bbb {R}^d)$.  
For the duration of this proof, we abbreviate $B(p, \rho_0/8)$ to $B$. There are two cases. 
\medbreak
Case one: for all  $\epsilon > 0$,  
\[
\omega(B\setminus \overline {W_{\epsilon}}, 
\partial (B\cap W_{\epsilon}), x) = 0
\]
for all $x\in B\setminus \overline {W_{\epsilon}}$. Then by Lemma~\ref{l:2.1}, $S$ is bounded
 from above and by Lemma~\ref{l:1.0}, $\tilde\Delta_{*} f_2<\infty$ everywhere. Also $B(p, p_0)$ was chosen so that $S$
is upper semi-continuous when restricted to $B(p,p_0)\cap \overline W$. Thus all the hypotheses of Lemma~\ref{l:2.2} are
satisfied and $W\cap B = \emptyset$, which is a  contradiction.
\medbreak
Case two: for some $\epsilon > 0$ and for some 
$x_0\in B\setminus \overline {W_{\epsilon}}$, we have
\begin{equation}\label{2.10}
\omega(B\setminus \overline {W_{\epsilon}}, 
\partial (B\cap W_{\epsilon}), x_0) > 0.
\end{equation}

Even though $\overline U$ is upper semi-continuous everywhere, it may
not be uniformly upper semi-continuous on $W_{\epsilon}$. This presents
a problem which did not arise at the corresponding point in Bourgain's
proof. To deal with this, we now introduce a subset of $W_{\epsilon}$
called $Z_{\epsilon}$, on a portion of which there holds a kind of
uniform upper semi-continuity.

Let 
\[
Z_\epsilon = \{y\in B\cap W_\epsilon: 
\omega(B\setminus \overline {W_{\epsilon}}, 
\overline{B(y, \delta)\cap W_{\epsilon}}, x_0)>0, \hbox{ for all }
\delta>0\}.
\]
Then, 
\begin{equation}\label{2.11}
\omega(B\setminus \overline {W_{\epsilon}}, 
\partial Z_\epsilon, x_0)>0.
\end{equation}
In fact, by definition, for each $z\in B\cap
W_\epsilon\setminus \overline {Z_\epsilon}$, there exists a ball 
$B(z, \delta_z)$, such that 
\begin{equation}\label{2.12}
\omega(B\setminus \overline {W_{\epsilon}}, 
\overline{B(z, \delta_z)\cap W_{\epsilon}}, x_0)=0. 
\end{equation}
The open cover $\{B(z, \delta_z)\}$ of $B\cap
W_\epsilon\setminus \overline {Z_\epsilon}$ has a countable subcover $\{B(z_i,\delta_{z_i})\}$. 
Thus, \eqref{2.12} implies
\[
\omega(B\setminus \overline {W_{\epsilon}}, 
B\cap W_\epsilon\setminus \overline {Z_\epsilon}, x_0)=0 .
\]
So \eqref{2.11} follows from \eqref{2.10} as $\omega(B\setminus \overline {W_{\epsilon}}, 
\overline {Z_\epsilon}, x_0)=\omega(B\setminus \overline {W_{\epsilon}}, 
\partial {Z_\epsilon}, x_0)$.

Since $\overline U$ is upper semi-continuous, 
\[
\bigcup_{m\geq 1}\left\{y\in \overline {Z_\epsilon}: \sup_{|z-y|\leq
2/m} \overline U(z) -\overline U(y) \leq  \frac{\epsilon}{40}\right\} =
\overline {Z_\epsilon}.
\]
Apply Baire's category theorem to the space $\overline {Z_\epsilon}$ to see that there exists $m\geq 1$ and an open ball 
$B(q, \rho)\subset B, q\in Z_\epsilon$, 
such that 
\[
B(q, \rho)\cap \overline {Z_\epsilon} \subset \overline{\left\{y\in 
\overline {Z_\epsilon}: \sup_{|z-y|\leq  2/m}
\overline U(z) -\overline U(y) \leq  \frac{\epsilon}{40}\right\}}.
\]
Equivalently, for any fixed $y\in B(q, \rho)\cap \overline {Z_\epsilon}$, 
there exists a sequence $y_n\in B(q, \rho)\cap Z_\epsilon$ convergent to $y$ such that  
\begin{equation}\label{2.13}
\sup_{|z-y_n|\leq  2/m} \overline U(z) - \overline U(y_n) \leq 
\frac{\epsilon}{40}.
\end{equation}
However, $\overline U$ is upper semi-continuous. So there exists $0< \delta <
\frac{1}{m}$ such that for $|y_n-y|< \delta$,  
\[
\overline U(y_n) - \overline U(y) < \frac{\epsilon}{40}.
\]
Thus, for $|z-y|< \frac{1}{m}$, since $|z-y_n|<\frac{2}{m}$ if
$|y_n-y|<\delta$, 
\begin{equation}\label{2.14}
\sup_{|z-y|<\frac{1}{m}}\overline U(z) < \overline U(y_n) +
\frac{\epsilon}{40} < \overline U(y) +\frac{\epsilon}{20}.
\end{equation}
Without loss of generality, we assume that $\frac{1}{m} \leq 
\frac{\rho}{2}$. 

Because $q\in Z_\epsilon$, we also have 
\begin{equation}\label{2.15}
\omega(B\setminus \overline {W_{\epsilon}}, 
\overline {B(q, \frac{\rho}{2})\cap W_\epsilon}, x_0)=\omega(B\setminus \overline {W_{\epsilon}}, 
\partial (B(q, \frac{\rho}{2})\cap W_\epsilon), x_0) >0.
\end{equation}

Set $F_\epsilon = B(q, \frac{\rho}{2})\cap W_\epsilon$. Then the rightmost inequality of 
Lemma~\ref{l:2.0} and \eqref{2.15} imply that 
\[
\omega(B\setminus \overline F_{\epsilon}, 
\partial F_\epsilon, x_0) \geq \omega(B\setminus \overline 
W_{\epsilon},  \partial F_\epsilon, x_0) >0.
\]

From Lemma~\ref{l:2.3}, there exists $p'\in \overline F_{\epsilon}$ 
such that 
\begin{equation}\label{2.16}
\inf_{\delta_1> 0}\underset{\delta_2\to 0}\liminf\inf_{x\in B(p',
\delta_2)} \, \omega(B(p', \delta_1)\setminus\overline F_{\epsilon}, 
\partial  (B(p',\delta_1)\cap F_{\epsilon}), x) = 1.
\end{equation}

Notice that Lemma~\ref{l:2.3} requires the set $F$ to be closed, so we cannot be sure that $p'\in F_\varepsilon$. 
Although $F_\varepsilon$ may not be closed, the uniformity implied by \eqref{2.16} allows us to continue.

Since $f_2$ restricted to $\overline W \cap B$ 
is continuous, we may select $1/(8m) > \delta_1>0$ such that 
\begin{equation}\label{2.17}
|f_2(z) - f_2(y)| \leq  \frac{\varepsilon}{10}
\end{equation}
for all $y, z\in B(p', 8\delta_1)\cap \overline W$.

Let $\eta>0$ be any positive number. From \eqref{2.16}, it follows that 
there exists $0< \delta_2 = \delta_2(\eta, \delta_1)< \delta_1$ such that 
\[
\omega(B(p', \delta_1)\setminus \overline F_{\epsilon}, 
\partial (B(p', \delta_1)\cap
F_{\epsilon}), y) > 1-\eta\hbox{ for all }y\in B(p', \delta_2).
\]

We may also assume that $\delta_1+\delta_2 = \delta_3^{'}$ satisfies
$B(p', \delta_3^{'})\subset B(p, \frac{\rho_0}{2})$. Pick any $\delta_3$ bigger than $\delta'_3$ but small enough to force $B(p', \delta_3) \subset B(p, \frac{\rho_0}{2})$. Note
that 
$p'\in \overline F_\varepsilon$ implies that there exists 
$p_1\in B(p', \frac{\delta_2}{2})\cap F_\epsilon$. Since 
$B(p', \delta_1) \subset B(p_1, \delta_3)$, 
\[
B(p_1, \delta_3)\setminus\overline{[F_\varepsilon\cup \{B(p_1, \delta_3)\setminus B(p', \delta_1)\}]}=B(p',\delta_1)\setminus\overline F_\varepsilon  .
\]
So by the rightmost inequality of Lemma~\ref{l:2.0},
\begin{eqnarray*}
&& \hskip-48pt \omega(B(p_1, \delta_3)\setminus \overline F_{\epsilon}, 
 \partial (B(p', \delta_1)\cap
F_{\epsilon}), y)\\ \noalign{\vskip4pt}
\quad & \geq& \omega (B(p_1, \delta_3)\setminus\overline{[F_\varepsilon\cup \{B(p_1, \delta_3)\setminus B(p',
\delta_1)\}]},
\partial (B(p', \delta_1)\cap F_{\epsilon}), y) \\ \noalign{\vskip4pt}
\quad & =& \omega(B(p', \delta_1)\setminus \overline F_{\epsilon}, 
\partial (B(p', \delta_1)\cap
F_{\epsilon}), y)\\ \noalign{\vskip4pt}
\quad &\geq& 1-\eta\hbox{ for all }y\in B(p_1, \frac{\delta_2}{2}).
\end{eqnarray*}
Consequently,
\[
\omega(B(p_1, \delta_3)\setminus \overline F_{\epsilon}, 
\partial (B(p_1, \delta_3)\cap
F_{\epsilon}), y) \geq  1-\eta\hbox{ for all }y\in B(p_1,
\frac{\delta_2}{2}),
\]
since $B(p', \delta_1)\cap F_{\epsilon} \subset B(p_1,
\delta_3)\cap F_{\epsilon}$.
Finally, by the left inequality of Lemma~\ref{l:2.0}
\[
\omega(B(p_1, \delta_3)\setminus \overline W, 
\partial (B(p_1, \delta_3)\cap
W), y) \geq \omega(B(p_1, \delta_3)\setminus \overline F_{\epsilon}, 
\partial (B(p_1, \delta_3)\cap F_{\epsilon}), y).
\]
We therefore have 
\begin{equation}\label{2.18}
\omega(B(p_1, \delta_3)\setminus \overline W, 
\partial (B(p_1, \delta_3)\cap
W), y)> 1-\eta\hbox{ for all }y\in B(p_1,
\frac{\delta_2}{2}).
\end{equation}
By definition, $p_1\in W_{\epsilon}$ implies that there exists $p_2\in
B(p_1, \frac{\delta_2}{2})$ such that 
\[
S(p_2)-S(p_1) \geq \frac{\varepsilon}{2}.
\]
Apply Lemma~\ref{l:2.1} at $p_1, p_2,$ and $\rho_1 = \delta_3$ where the 
inequality \eqref{2.4} holds for $\delta_3$. Then by \eqref{2.14},  
\eqref{2.17}, \eqref{2.18}, and the above inequality, we 
have 
\begin{equation}\label{2.19}
\frac{\varepsilon}{2} \leq   S(p_2) - S(p_1) \leq  c[|f_2(p_1)| +
\delta_3^{-\frac34 (d-1)}]\eta^{1/4} + \frac{\varepsilon}{5}.
\end{equation}
Note here that $p_1, p_2$, and $\delta_3$ depend on $\eta$. However, since $f_2$ is continuous and hence bounded on $B(p,
\rho_0)\cap\overline W$  and $\delta_3$ is bounded below by $\delta_1$ as $\eta \to 0$, so \eqref{2.19} becomes a
contradiction upon choosing $\eta$ sufficiently small.

\section{Proof of Lemma~\ref{l:2.1}}\label{s:3}

For any  bounded measurable function $f(x)$ defined on $\partial G$, 
\begin{equation}\label{3.1}
H_f(x) = \int_{\partial G} f(z)\, \omega(G, dz, x) 
\end{equation}
is harmonic in $G$. If every point on $\partial G$ satisfies the  exterior 
cone condition and $f$ is continuous at $x\in \partial G$, then 
\[
\lim_{{y\rightarrow x\atop y\in G}} H_f(y)=f(x).
\] 
Since any upper semi-continuous function is the limit of a decreasing sequence of continuous functions, so the maximum principle for subharmonic functions and \eqref{3.1} imply that 
\begin{equation}\label{3.1.1}
f(x) \leq \int_{\partial G} f(z)\, \omega(G, dz, x) 
\end{equation}
for any function $f$ subharmonic on an open set $\tilde G\supset \overline G \supset G$.

We need only to consider $p_2\notin \overline W$. Let $\tau = \hbox{dist}(p_2, \overline W) \leq  \frac12 \rho_1$. For $\kappa \ll \tau$, define
\[
G_\kappa = \{x\in B(p_1, \rho_1): \hbox{dist} (x, \overline W) < \kappa\}.
\]
Clearly $\overline W \cap B(p_1, \rho_1)\subset G_\kappa$. We know that 
$S$ is upper semi-continuous and $\tilde \Delta^* S(x) \geq 0$ on 
$B(p, \rho_0)\setminus\overline W$. This is the hypothesis of a 
classical theorem (see for example \cite[p.~14]{R}) which concludes that 
$S$ is subharmonic on $B(p, \rho_0)\setminus \overline W$. Thus, 
$S(x) - S(p_1)$ is subharmonic on $B(p, \rho_0)\setminus \overline
W$. In particular, $S(x) - S(p_1)$ is subharmonic on an open set containing $\overline {B(p_1, \rho_1)
\setminus \overline G_\kappa}$. Note that $B(p_1, \rho_1)\setminus \overline
G_\kappa$ satisfies the exterior cone condition everywhere on the boundary.
So by \eqref{3.1.1}, we 
have 
\begin{eqnarray*}
\noalign{\vskip6pt}
&&\hskip-.5in S(p_2) - S(p_1)   \leq    \int\limits_{\partial (B(p_1,
\rho_1)\setminus \overline G_\kappa)}\, 
[S(x) - S(p_1)]\, \omega(B(p_1, \rho_1)\setminus
\overline G_\kappa, dx, p_2)\\ \noalign{\vskip6pt}
&&\quad =\ \int\limits_{\partial B(p_1, \rho_1)\setminus (B(p_1, \rho_1)\cap 
\partial G_{\kappa})}\,
[S(x) - S(p_1)]\, \omega(B(p_1, \rho_1)\setminus
\overline G_\kappa, dx, p_2)\\ \noalign{\vskip6pt}
&&\qquad+ \ \int\limits_{B(p_1, \rho_1)\cap\partial G_\kappa}\,
[S(x)-S(p_1)]\, \omega(B(p_1, \rho_1)\setminus
\overline G_\kappa, dx, p_2)\\ \noalign{\vskip6pt}
&&\quad  =\ I_1 + I_2.
\end{eqnarray*}

We first estimate $I_1$.  When $p_2\in B(p_1, \rho_1/2)$, a classical result on
harmonic measure shows that $\omega = 
\omega(B(p_1, \rho_1)\setminus \overline G_\kappa, dx, p_2)$ is absolutely
continuous with respect to the surface Lebesgue measure $\sigma$ when 
restricted to the sphere $B(p_1, \rho_1)$. (See \cite{D} or
(4.39) of \cite{AW}.) By \eqref{1.2}, $f_2(x) = F(x)$ 
almost everywhere with respect to Lebesgue measure; thus, for almost every 
$\rho_1>0$, $f_2(x) = F(x)$ almost everywhere with respect to the surface Lebesgue measure on $B(p_1, \rho_1)$ and hence with respect to the harmonic measure $\omega$ for all $p_2\in B(p_1, \rho_1/2)$. 
Consequently,
\begin{eqnarray*}
\noalign{\vskip4pt}
I_1 &\leq&  \int\limits_{\partial B(p_1, \rho_1)} |F(p_2) - f_2(p_1)|\,
\omega(B(p_1, \rho_1)\setminus \overline G_\kappa, dx, p_2)\\ \noalign{\vskip4pt}
&&+\ \underset{q\in B(p_1, \rho_1)}\sup \overline U(q) -\overline
U(p_1)\\ \noalign{\vskip4pt} 
& =& I_3 + I_4. 
\end{eqnarray*}
A result of Bourgain \cite{B}  (see also Lemma 4.5 of \cite{AW})  shows that 
\[
I_3 \leq  c\left( \left[|f_2(p_1)|
+\rho_1^{-\frac34(d-1)}\right] [1-\omega (B(p_1, \rho_1)\setminus\overline
W, \partial(W\cap B(p_1, \rho_1)), p_2)]^{\frac14}\right).
\]
\demo{{R}emark\/ {\rm 4.1}} 
In fact, Lemma 4.5 of \cite{AW} was based on Connes' condition
\eqref{0.3}. But a careful reading of the proof shows the conclusion of
  Lemma 4.5 holds true under Bourgain's condition \eqref{e:bou} since
inequalities (4.18) and (4.19) in Lemma 4.2 and Corollary 4.3
respectively can be replaced by 
\[
\sup_{k} \frac{1}{2^{k}}\sum_{|\xi|^2\sim 2^k}|c_{\xi}|^2 \leq M
\]
as they are used only in (4.21).
\enddemo 

This gives the first half of \eqref{2.4}.  Now we estimate $I_2$. 

For any $x\in B(p_1, \rho_1)\cap \partial G_{\kappa}$, there exists
$\tilde x \in \overline W\cap B(p_1, 2\rho_1)$, such that $|x-\tilde x| =
\kappa$. Since $S$ is subharmonic at $x$, 
\begin{equation}\label{3.2}
S(x) \leq  A_\kappa f_2(x) + A_\kappa \overline U(x).
\end{equation}
Since $\tilde x\in \overline W\cap B(p, \rho_0)$, by assumption,  
\begin{equation}\label{3.3}
|f_2(\tilde x) - f_2(\tilde x, \kappa)| \leq N \kappa.
\end{equation}

Thus combining \eqref{3.2} and \eqref{3.3}, we have 
\begin{eqnarray*}
S(x) - S(\tilde x) &\leq&  A_{\kappa} f_2(x) - 
A_{\kappa} f_2(\tilde x, \kappa) + A_{\kappa} f_2(\tilde x, \kappa)-
f_2(\tilde x, \kappa)\\
&& +\ A_{\kappa} \overline U(x)-\overline U(\tilde x) + N\kappa.
\end{eqnarray*}
Consequently, 
\begin{eqnarray*}
S(x) - S(p_1) & = &  S(x) - S(\tilde x) + S(\tilde x) - S(p_1)\\ \noalign{\vskip4pt}
&\leq&  A_\kappa f_2(x) - A_\kappa f_2(\tilde x, \kappa) + 
A_{\kappa} f_2(\tilde x, \kappa)- f_2(\tilde x, \kappa)\\ \noalign{\vskip4pt}
&& +\  f_2(\tilde x) -f_2(p_1) + A_\kappa\overline U(x) - \overline
U(p_1) + N\kappa\\ \noalign{\vskip4pt}
&\leq&  A_\kappa f_2(x) - A_\kappa f_2(\tilde x, \kappa) + A_{\kappa}
f_2(\tilde x, \kappa)- f_2(\tilde x, \kappa)\\ \noalign{\vskip4pt}
&& +\  \underset{q\in B(p_1, 2\rho_1)\cap \overline W}\sup |f_2(q)-
f_2(p_1)|\\ \noalign{\vskip4pt}
&& +\ \underset{q\in B(p_1, 2\rho_1)}\sup \overline U(q)- 
\overline U(p_1) + N\kappa.
\end{eqnarray*} 
From \eqref{1.2} and the definition of $A_\rho f_2(x)$, we have $A_\rho
f_2(x) = A_\rho F(x)$ for all $x$. Thus 
\begin{eqnarray}\label{3.4}
&& \\ 
I_2 &\leq& \int\limits_{B(p_1, \rho_1)\cap \partial G_\kappa} 
|A_\kappa F(x) - A_\kappa f_2(\tilde x, \kappa)| \, 
\omega(B(p_1, \rho_1)\setminus \overline G_\kappa,
dx, p_2)\nonumber \\ \noalign{\vskip4pt}
&&+ \ \int\limits_{B(p_1, \rho_1)\cap \partial G_\kappa} 
|A_\kappa f_2(\tilde x, \kappa) - f_2(\tilde x, \kappa)| \, 
\omega(B(p_1, \rho_1)\setminus \overline G_\kappa,
dx, p_2)\nonumber \\ \noalign{\vskip4pt}
&& +\  \underset{q\in B(p_1,
2\rho_1)\cap \overline W} \sup |f_2(q)-f_2(p_1)|\nonumber \\
&& +\ \underset{q\in B(p_1, 2\rho_1)}\sup
\overline U(q)- \overline U(p_1)  + N\kappa\nonumber\\ \noalign{\vskip4pt}
& = &I_5 + I_6 + \underset{q\in B(p_1,
2\rho_1)\cap \overline W} \sup |f_2(q)-f_2(p_1)|\nonumber \\ \noalign{\vskip4pt}
&&+\  \underset{q\in B(p_1, 2\rho_1)}\sup
\overline U(q)- \overline U(p_1) + N\kappa. \nonumber
\end{eqnarray} 

It is enough to show that $I_5\to 0$ and $I_6\to 0 $ as $\kappa\to 0$. 
Observe that 
\begin{eqnarray*}
I_5 &\leq& \int\limits_{B(p_1, \rho_1)\cap \partial G_\kappa} 
|A_\kappa F(x) - A_\kappa F(\tilde x)| \, 
\omega(B(p_1, \rho_1)\setminus \overline G_\kappa,
dx, p_2)\\
&& + \ \int\limits_{B(p_1, \rho_1)\cap \partial G_\kappa} 
|A_\kappa F(\tilde x) - A_\kappa f_2(\tilde x, \kappa)| \, 
\omega(B(p_1, \rho_1)\setminus \overline G_\kappa,
dx, p_2)\\
& =& I_7 +I_8.
\end{eqnarray*}
A result of Bourgain \cite{B} (again, see also Lemma 4.4 of \cite{AW}),
shows that $I_7\to 0$ as $\kappa\to 0$. However, the same proof of 
Lemma 4.4 in \cite{AW} shows that if $a_{\xi}$ satisfies Bourgain's
condition \eqref{e:bou}, then 
\begin{equation}\label{3.5}
\lim_{\kappa\to 0} \kappa \int\limits_{B(p_1, \rho_1)\cap \partial G_\kappa}
\left|\sum_{|\xi|\neq 0} \frac{a_{\xi}}{|\xi|}\hat I(\kappa |\xi|) 
e^{i\langle \overline x, \xi\rangle}\right| \, 
\omega(B(p_1, \rho_1)\setminus \overline G_\kappa, dx, p_2) = 0,
\end{equation}
where $|\overline x - x|\leq \kappa$.
Note that
\[
A_\kappa F(\tilde x) - A_\kappa f_2(\tilde x, \kappa) = -\sum_{|\xi|\neq 0} 
\frac{a_{\xi}}{|\xi|^2}\hat I(\kappa |\xi|)e^{i\langle\tilde x, \xi\rangle}
(1-e^{-|\xi|\kappa}),
\]
while by the mean value theorem, for each $\xi\neq
0$, there exists $t_{\xi} > 0$ such that 
\[
\sum_{|\xi|\neq 0} \frac{a_{\xi}}{|\xi|^2}\hat I(\kappa |\xi|)e^{i\langle 
\tilde x, \xi\rangle}(1-e^{-|\xi|\kappa}) = \kappa \sum_{|\xi|\neq 0}\frac{a_{\xi}
e^{-|\xi|t_{\xi}}}{|\xi|}\hat I(\kappa |\xi|) 
e^{i\langle\tilde x, \xi\rangle}.
\]
Since $e^{-|\xi|t_{\xi}} < 1$, $\{a_{\xi}e^{-|\xi|t_{\xi}}\}$ satisfies
Bourgain's condition \eqref{e:bou} as $\{a_{\xi}\}$ does. Thus by
\eqref{3.5}, $I_8\to 0$ as $\kappa\to 0$. This shows that $I_5\to 0$ as
$\kappa \to 0$. 

The method that Bourgain used to prove that $I_7\to 0$ as $\kappa\to 0$
can also be used to prove $I_6\to 0$ as $\kappa\to0$. To establish this, we will use the following lemma of Bourgain \cite {B}. (See also the proof of Corollary~4.3 of \cite{AW}.)
\proclaim{Lemma}\label{l:3.1}
Let $k\geq 1, \gamma>0, \eta\leq 2^{-k}${\rm .} Let $E_{k, \gamma, \eta}$ be
a set of $\eta$\/{\rm -}\/separated points $x\in B(p, q)\subset \Bbb  R^d$
satisfying
\[
\left|\sum_{|\xi|\sim 2^k} \frac{b_{\xi}}{|\xi|^2} e^{i\langle x,
\xi\rangle}\right|\geq \gamma.
\]
Then{\rm ,} the cardinality of $E_{k,\gamma, \eta}$ satisfies
\[
|E_{k,\gamma,\eta}|\leq c \gamma^{-2}\eta^{-d}2^{-2k}\nu_k^2,
\]
where $c$ is an absolute constant and 
\[
\nu_k^2 = 2^{-2k}\sum_{|\xi|\sim 2^{k}} |b_{\xi}|^2.
\]
\endproclaim
Let
\[
\alpha_{k} = \left\{ \begin{array}{ll} c\alpha [\log(1+2^k \kappa)]^{-2}, &
\hbox{ for } 2^k\geq \kappa^{-1}\\
                           c\alpha [\log(1+2^{-k} \kappa^{-1})]^{-2}, &
\hbox{ for } 2^k < \kappa^{-1} .\end{array} \right. 
\]
The positive constant $c$ is chosen so that $\sum_{k\geq 1}\alpha_{k}
\leq 2c\alpha\sum_{n\geq 0}(\log(1+2^n))^{-2}\break = \alpha$ for all $\alpha> 0$. Clearly, $c$ is an absolute constant. 
For $\alpha >0$, let
\[
S_{\kappa, \alpha} = \{ x\in B(p_1, \rho_1)\cap \partial G_{\kappa}: 
|A_{\kappa}f_2(\tilde x, \kappa) - f_2(\tilde x, \kappa)| > \alpha\}.
\]
Then
\begin{equation}\label{3.6}
I_6 = \int_{0}^{\infty}\omega(B(p_1, \rho_1)\setminus\overline
G_{\kappa}, S_{\kappa, \alpha}, p_2)\, d\alpha.
\end{equation}

Let
\begin{eqnarray*}
S_{\kappa, k,\alpha_{k}} & = &  \{x\in B(p_1, \rho_1)\cap \partial G_{\kappa}: 
|A_{\kappa}f_{2, k}(\tilde x, \kappa) - f_{2, k}(\tilde x, \kappa)| > 
\alpha_{k}\},\\ \noalign{\vskip4pt}
S'_{\kappa, k, \alpha_{k}} & =& \{x\in B(p_1, \rho_1): 
|A_{\kappa}f_{2, k}(x, \kappa) - f_{2, k}(x, \kappa)| > 
\alpha_{k}\},
\end{eqnarray*}
where
\begin{eqnarray*}
f_{2, k}(x, \kappa) & = &  \sum_{|\xi|\sim 2^k} \frac{a_{\xi}}{|\xi|^2}
e^{i\langle x, \xi\rangle-\kappa|\xi|}, \hbox{ and} \\ \noalign{\vskip4pt}
A_{\kappa}f_{2, k}(x, \kappa)& = & \sum_{|\xi|\sim 2^k} \frac{a_{\xi}}{|\xi|^2}\hat I(\kappa|\xi|)
e^{i\langle x, \xi\rangle-\kappa|\xi|}.
\end{eqnarray*}
Then
\begin{equation}\label{3.7}
S_{\kappa, \alpha} \subset \bigcup_{k\geq 1} S_{\kappa, k, \alpha_k}.
\end{equation}
Since $|x-\tilde x|= \kappa$, observe that a collection of balls of radius 
$\eta\leq 2^{-k}$ centered at points in $S'_{\kappa, k, \alpha_k}$ covering  
$S'_{\kappa, k, \alpha_k}$ will cover $S_{\kappa, k, \alpha_k}$ if the
radius of each ball is enlarged by $\kappa$. 

Bourgain's condition \eqref{e:bou} may be restated as $\delta_k\to 0$, where 
\[
\delta^2_k=2^{-2k}\sum_{|\xi|\sim 2^k}|a^2_{\xi}|.
\]
In particular, $\delta^2_k$ is bounded for all $k$. Now apply Lemma~\ref{l:3.1} and use the fact that $\hat I(|\xi|) = O(|\xi|^{-(d+1)/2})$ as $|\xi|\to \infty$ and 
also use \eqref{1.3}. We find that the number of balls of radius $\eta\leq 2^{-k}$ centered at $S'_{\kappa, k, \alpha_k}$
covering
$S'_{\kappa, k,\alpha_k}$ is at most
\begin{eqnarray}\label{3.8}
|E_{k, \alpha_k, \eta}| 
&\leq &c\alpha_k^{-2}\eta^{-d}2^{-2k}  e^{-\kappa 2^k}\delta_k^2\sup_{2^{k-1}\leq j <2^{k}}|\hat I(\kappa j) -1|^2\\
& \leq & \left\{
\begin{array}{ll} c\alpha_k^{-2} \delta_k^2 \eta^{-d} 2^{-2k} e^{-
\kappa 2^k}, & \quad  \kappa 2^k \geq 1\nonumber \\ \noalign{\vskip4pt}
c\alpha_k^{-2} \delta_k^2 \eta^{-d} 2^{-2k} (\kappa 2^{k})^4,  & \quad
\kappa 2^k < 1. \end{array} \right. 
\end{eqnarray}

We estimate $\omega(B(p_1, \rho_1)\setminus \overline G_{\kappa}, 
S_{\kappa, k,\alpha_k}, p_2)$ according to the size of $k$.  

\demo{Case {\rm (i):} $\kappa 2^k \geq 1$}
By \eqref{3.8} with $\eta= 2^{-k}$ and the observation made after 
\eqref{3.7},  the number of balls of radius $2\kappa$ covering $S_{\kappa, k, \alpha_k}$ is at most
 $M = c\alpha_k^{-2} \delta_k^2 2^{(d-2)k}e^{-\kappa 2^{k}}$. Let $\{B_i\}_{1\le i\le M_1}, M_1\le M$,  denote these balls.
Then 
\begin{eqnarray}\label{3.9}
&& \nonumber \\
\noalign{\vskip-24pt}
&&\\
\omega(B(p_1, \rho_1)\setminus \overline G_{\kappa}, S_{\kappa, k,
\alpha_k}, p_2) &\leq &\omega(B(p_1, \rho_1)\setminus S_{\kappa, k,
\alpha_k}, S_{\kappa, k,
\alpha_k}, p_2)\nonumber \\ \noalign{\vskip5pt}
&\leq&\omega(B(p_1, \rho_1)\setminus \overline{\cup B_i}, \partial(\cup B_i), p_2)\nonumber \\ \noalign{\vskip5pt}
& = & \sum_{i=1}^{M_1} \omega(B(p_1, \rho_1)\setminus \overline{\cup B_i}, \partial{B_i}, p_2)\nonumber \\ \noalign{\vskip5pt}
&\leq& \sum_{i=1}^{M_1} \omega(B(p_1, \rho_1)\setminus \overline{B_i}, \partial{B_i}, p_2)\nonumber \\
\noalign{\vskip5pt} 
&\leq& c\alpha_k^{-2} \delta_k^2 2^{(d-2)k}e^{-\kappa
2^{k}}\left(\frac{\kappa}{\tau}\right)^{d-2}\nonumber \\ \noalign{\vskip5pt}
& =& \frac{c}{\tau^{d-2}}\alpha_k^{-2} \delta_k^2 (\kappa 2^{k})^{d-2} 
e^{-\kappa 2^k}\nonumber \\ \noalign{\vskip5pt}
& \leq& \frac{c}{\tau^{d-2}}\alpha_k^{-2} \delta_k^2 (\kappa 2^{k})^{-1}, \nonumber
\end{eqnarray} 
where $c$ is a constant which may vary from line to line. The first line
of \eqref{3.9} follows from the rightmost inequality of Lemma~\ref{l:2.0},
since $B(p_1, \rho_1)\setminus \overline G_{\kappa}$ has been relaced by
$B(p_1, \rho_1)\setminus S_{\kappa, k,
\alpha_k}$. The second line follows from the left-most inequality of
Lemma~\ref{l:2.0}, since both occurrences of $S_{\kappa, k, \alpha_k}$ have
been replaced by the union of balls $\bigcup_{i=1}^{M_1} B_i$ of radius
$2\kappa$ covering it. The third line follows from the subadditivity of
harmonic measure in the second coordinate. The fourth line follows from the
right-most inequality of Lemma~\ref{l:2.0}, since $\bigcup_{i=1}^{M_1} B_i$
has been replaced by $B_i$ in each term. To see the next line, write $B_i$
as $B(q_i, \rho_i)$; use the explicit formula for the Poisson integral to
estimate each term $\omega(B(p_1, \rho_1)\setminus \overline{B_i},
\partial{B_i}, p_2)$ by $\frac{\rho_i^{d-2}}{|p_2-q_i|^{d-2}}$, where
$\rho_i=2\kappa$ and $|p_2 -q_i|\geq \tau-3\kappa>\tau/2$; and finally use
the first line of \eqref{3.8} to estimate the number of terms. 
\enddemo

\vglue6pt 
\demo{Case {\rm (ii):} $\kappa^{-1/2}\leq 2^{k} < \kappa^{-1}$} 
By \eqref{3.8} with $\eta = 2^{-k}$, the number of balls of radius $2\cdot 2^{-k}$ covering $S_{\kappa, k, \alpha_k}$ is  at
most $c\alpha_k^{-2} \delta_k^2 2^{(d-2)k} (\kappa 2^{k})^{4}$. So as shown in Case (i),
\begin{eqnarray}\label{3.10}
&& \nonumber \\
\noalign{\vskip-24pt}
&&\\
\noalign{\vskip-6pt}
\omega(B(p_1, \rho_1)\setminus \overline G_{\kappa}, S_{\kappa, k,
\alpha_k}, p_2) &\leq& c\alpha_k^{-2} \delta_k^2 2^{(d-2)k}(\kappa
2^{k})^4\left(\frac{1}{2^{k}\tau}\right)^{d-2}\nonumber \\ \noalign{\vskip5pt}
& \leq&\frac{c}{\tau^{d-2}}\alpha_k^{-2} \delta_k^2 (\kappa 2^k)^{2}. \nonumber
\end{eqnarray}
\enddemo
\demo{Case {\rm (iii):} $2^{k}< \kappa^{-1/2}$}
By \eqref{3.8} with $\eta = \sqrt{\kappa}$, the number of balls of radius $2\sqrt{\kappa}$ covering $S_{\kappa, k,
\alpha_k}$ is at most $\alpha_k^{-2} \delta_k^2 2^{-2k} \kappa^{-d/2}(\kappa 2^{k})^{4}$. So, 
\begin{eqnarray}\label{3.11}
&&\\ 
\noalign{\vskip-6pt}
\omega(B(p_1, \rho_1)\setminus \overline G_{\kappa}, S_{\kappa, k,
\alpha_k}, p_2) &\leq& c\alpha_k^{-2} \delta_k^2 2^{-2k}\kappa^{-d/2}(\kappa
2^{k})^4\left(\frac{\sqrt{\kappa}}{\tau}\right)^{d-2}\nonumber \\ \noalign{\vskip4pt}
&\leq& c\alpha_k^{-2} \delta_k^2 2^{-2k}\kappa^{-d/2}(\kappa
2^{k})^2\left(\frac{\sqrt{\kappa}}{\tau}\right)^{d-2}\nonumber \\ \noalign{\vskip4pt}
&\leq& \frac{c}{\tau^{d-2}}\alpha_k^{-2} \delta_k^2 \kappa. \nonumber
\end{eqnarray}
When \eqref{3.9}, \eqref{3.10}, and \eqref{3.11} are combined, it follows from
\eqref{3.7} and the definitions of $\alpha_k$ that
\begin{eqnarray}\label{3.12}
&&\\ 
\noalign{\vskip-6pt}
 \omega(B(p_1, \rho_1)\setminus \overline G_{\kappa}, S_{\kappa, \alpha},
p_2)   &\leq &\frac{c}{\tau^{d-2}}\left\{\sum_{2^k \geq \kappa^{-1}}
\alpha_k^{-2}\delta_k^2(2^k\kappa)^{-1} +  \right. \nonumber\\
\noalign{\vskip4pt}
&& 
\sum_{\kappa^{-1/2}\leq 2^k <\kappa^{-1}} \left.\alpha_k^{- 2}\delta_k^2(2^k\kappa)^{2}+
\sum_{2^k<\kappa^{-1/2}}\alpha_k^{- 2}\delta_k^2\kappa\right\}\nonumber \\ \noalign{\vskip4pt}
&\leq& \frac{c}{\tau^{d-2}}\alpha^{-2} \left\{\max_{2^k>\kappa^{-
1/2}}\delta_k^2 + \kappa |\log\kappa|^5\right\}. \nonumber
\end{eqnarray}
Choose 
\[
\beta_{\kappa} = \left\{\max_{2^k>\kappa^{-1/2}}\delta_k^2 + 
\kappa |\log\kappa|^5\right\}^{1/2}.
\]
Then by \eqref{3.6}, 
\[
I_6 \leq \int_{0}^{\beta_\kappa}\, d\alpha + \frac{c}{\tau^{d-
2}}\beta_\kappa^2\int_{\beta_\kappa}^{\infty} \alpha^{-2}\, d\alpha
 = \left(\frac{c}{\tau^{d-2}}+1\right)\beta_{\kappa}\to 0 
\]
as $\kappa\to 0$. This completes the proof. 

\section{Proof of Lemma~\ref{l:2.2}}\label{s:4}

Let the average of $H$ on the surface of $B(x, \rho)$  be denoted by   
\[
D_\rho H(x) = \frac{1}{\sigma_d \rho^{d-1}} \int_{\partial B(x,
\rho)} H(z)\, d\sigma(z),
\]
where $\sigma$ is the surface measure, and $\sigma_d = d v_d$ is the surface 
area of the unit ball in $\Bbb {R}^d$. Then, 

\begin{eqnarray}
&&\nonumber \\
\noalign{\vskip-32pt}
A_\rho H(x) & = &  \frac{1}{v_d \rho^d}\int_{B(x,\rho)} H(z)\, dz\label{4.1}\\
\noalign{\vskip4pt}
& = &\frac{1}{v_d\rho^d}\int_{0}^{\rho} \int_{\partial B(x, \beta)} H(z)\,
d\sigma(z)\, d\beta\nonumber \\ \noalign{\vskip4pt}
& =& \frac{\sigma_d}{v_d\rho^d}\int_{0}^{\rho} D_\beta H(x)\beta^{d-1}\,
d\beta \nonumber \\ \noalign{\vskip4pt}
& = &\frac{d}{\rho^d}\int_{0}^{\rho} D_\beta H(x)\beta^{d-1}\,
d\beta .  \nonumber
\end{eqnarray} 
For any $\eta >0$ and $x\in B(p, r)$, by \eqref{2.5}  there
exist  two sequences $\rho_{i,n} = \rho_{i, x, \eta, n} \downarrow 0$ 
(with $\rho_{i,n} < r - |x-p|$), $i=1, 2$, such that for all $n\geq 1$,  
\[
A_{\rho_{1,n}} S(x)- S(x) \geq -\eta \rho_{1,n}^2 \hbox{ and } A_{\rho_{2,n}}f_2(x) - f_2(x) \leq c_x \rho_{2,n}^2,
\]
where $c_x$ is a positive constant independent of $\rho_{2,n}$.
Thus, by \eqref{4.1}, the above inequalities imply for all $n \geq 1$, 
\begin{eqnarray}
\int_{0}^{\rho_{1,n}}[ D_\beta S(x) - S(x) + a\eta \beta^2] \beta^{d-1}\, d\beta
&\geq& 0\\ \noalign{\vskip4pt}
\int_{0}^{\rho_{2,n}}[ D_\beta f_2(x) - f_2(x) - ac_x \beta^2] \beta^{d-1}\, d\beta
&\leq&  0 \nonumber
\end{eqnarray}
where $a=\frac{d+2}{d}$. So there exist $\beta_n =
\beta_{x, \eta, n}\leq  \rho_{1,n}, \beta_n \downarrow 0$, and 
$r_n = r_{x, \eta, n}\break\leq  \rho_{2,n}, r_n \downarrow 0$, such that 
\begin{equation}\label{4.2}
D_{\beta_n}S(x) - S(x) \geq -a\eta \beta_n^2\hbox{ and } D_{r_n}f_2(x) - f_2(x) \leq ac_x  r_n^2.
\end{equation}
Let $B(q, \rho_1) \subset B(p, r)$. We show, for $y\in B(q,
\rho_1)\setminus\overline W$, that  
\begin{eqnarray}\label{4.4}
S(y) & \leq&  \int_{\partial B(q, \rho_1)} S(z) \omega(B(q, \rho_1), dz, y)\\ \noalign{\vskip4pt}
& = &\frac{1}{\sigma_d \rho_1} \int_{\partial B(q, \rho_1)}
\frac{\rho_1^2-|q-y|^2}{|z-y|^d} S(z)\, d\sigma (z). \nonumber
\end{eqnarray}

If \eqref{4.4} holds, then for $q\in W$, by \eqref{4.2}, there exists a
decreasing sequence $r_n$ of positive numbers going to $0$ such that for
each $n$  
\begin{equation}\label{4.5}
f_2(q) - D_{r_n}f_2(q) \geq - ac_q r_n^2. 
\end{equation}
For any given $\epsilon>0$, using upper semi-continuity of $\overline U$ at 
$q$, we have for large~$n$,  
\[
\overline U(q) \geq \sup_{|y-q|\leq  r_n}\overline U(y) - \epsilon. 
\]
Thus, for large $n$, 
\[
\overline U(q) \geq D_{r_n}\overline U(q) -\epsilon.
\]
Consequently
\begin{equation}\label{4.6}
S(q) \geq D_{r_n} S(q)-ac_q r_n^2 -\epsilon.
\end{equation}
Note that for each $r>0$, by the mean value theorem, there exists a 
constant~$c$ such that
\[
\left|\frac{1}{r^{d-2}} -\frac{r^2 - |q-y|^2}{|z-y|^{d}}\right|\leq 
c\frac{|q - y|}{r^{d-1}},
\] 
if $|y-q|< \frac12 |z-q|=\frac12r$. Therefore, for $|y-q|<\frac12r$,  
\begin{eqnarray}\label{4.7}
&& \hskip-48pt \left|\frac{1}{\sigma_d r} \int_{\partial B(q, r)}
\frac{r^2-|q-y|^2}{|z-y|^d} S(z)\, d\sigma (z) - D_r S(q)\right|\\ \noalign{\vskip4pt}
&\leq &
\frac{1}{\sigma_d r}\int_{\partial B(q, r)}
\left|\frac{1}{r^{d-2}} -\frac{r^2 - |q-y|^2}{|z-y|^{d}}\right|\,
|S(z)| d\sigma(z)\nonumber \\ \noalign{\vskip4pt} 
&\leq&  c\frac{|q-y|}{r} D_r|S(q)|. \nonumber
\end{eqnarray} 
Combining \eqref{4.4}--\eqref{4.7}, we have for any given $\epsilon$, for
$n$ large,   
\[
S(y) - S(q)\leq  a c_q r_n^2 +c\frac{|q-y|}{r_n} D_{r_n} |S(q)| + \epsilon,
\]
if $|y-q|\leq  \frac12 r_n$ and $y\in B(q, r_n)\setminus \overline W$. 
Letting $y\to q$, then $n \to\infty$, and then $\epsilon \to 0$, we have
\[
\limsup_{{y\to q\atop  y \in B(q, \rho_1)\setminus\overline W}}
S(y) \leq  S(q).
\]
Thus, $S$ is upper semi-continuous at $q$ since $S$ is upper
semi-continuous when restricted to $B(p, r)\cap \overline W$.
Consequently, $W$ must be the empty set. So $S$ is upper semi-continuous 
in $B(p, r)$. Inequality \eqref{4.4} also implies that $S(q) \leq  A_\rho
S(q)$ for all $B(q, \rho) \subset B(p, r)$. Thus $S$ is subharmonic in
$B(p, r)$ since it is also in $L^1$. 

It only remains to prove \eqref{4.4}.  Let $\{X_t\}_{t\geq 0}$ be the 
standard Brownian motion starting from a fixed point $y\in B(q, \rho_1)\setminus\overline W$ 
in the probability space $(\Omega, {\cal F}, P^{y})$. Define
\[
T = \inf\{t\geq 0: X_t\in \partial B(q, \rho_1)\}
\]
to be the exit time of $X_t$ from $B(q, \rho_1)$. Then by \eqref{3.1}, 
inequality \eqref{4.4} is equivalent to
\begin{equation}\label{4.8}
S(y) \leq  E^{y} [S(X_T)].
\end{equation}

We first show that for any stopping time $S\leq  T$,
\begin{equation}\label{4.9}
P^y(X_S\in  W) = 0.
\end{equation}
This is implied by 
\begin{equation}\label{4.91}
P^y(X_S\in \overline W_{\varepsilon})=0
\end{equation}
as $W\subset \bigcup\limits_{\varepsilon>0} W_{\varepsilon}$.
Let $R$ be the hitting time of $X_t$ with $\partial (B(q,
\rho_1)\setminus\overline W_{\varepsilon})$: 
\[
R = \inf\{t\geq 0: X_t \in \partial (B(q, \rho_1)\setminus\overline W_{\varepsilon})\}.
\]
Then $R \leq  T$. Since $y\in B(q, \rho_1)\setminus\overline W\subset B(q, \rho_1)\setminus\overline W_{\varepsilon}$
and, by assumption,
\begin{eqnarray}\label{4.11}
0 & = &  \omega(B(q, \rho_1)\setminus\overline W_{\varepsilon},
\partial (B(q, \rho_1)\cap\overline W_{\varepsilon}), y)\\ \noalign{\vskip4pt}
& =&P^y(X_R\in \partial (B(q, \rho_1)\cap\overline W_{\varepsilon}))\nonumber \\ \noalign{\vskip4pt}
& =& P^y(X_R\in\overline W_{\varepsilon}), \nonumber
\end{eqnarray}
we see that 
\begin{equation}\label{4.10}
P^y(X_R\in\overline W_{\varepsilon}) = 0.
\end{equation} 
Next, by definition,
\[
\{R<T\} \subset \{X_R\in \partial W_{\varepsilon}\}.
\]
So 
\begin{equation}\label{4.12}
P^y(R <T) \leq  P^y(X_R\in \partial W_{\varepsilon}) = 0.
\end{equation}
Thus, by \eqref{4.11} and \eqref{4.12} we have
\begin{eqnarray}\label{4.13}
\qquad \enspace P^y(X_T\in\overline W_{\varepsilon}) & = &  P^y(X_T \in\overline W_{\varepsilon}, R =T) + P^y
(X_T\in\overline W_{\varepsilon}, R < T)\\  \noalign{\vskip4pt}
& \leq&  P^y(X_R \in\overline W_{\varepsilon} ) + P^y (R< T)\nonumber \\ \noalign{\vskip4pt}
& =& 0. \nonumber
\end{eqnarray}
To show \eqref{4.91} for a general stopping time $S$, note that for any 
$\tau>0$, there exists an open set $G$ such that $\overline {W_{\varepsilon}} \subset G$ 
and
\begin{equation}\label{4.14}
\omega(B(q, \rho_1)\setminus\overline G, \partial (B(q, \rho_1)\cap G), y)
<\tau. 
\end{equation}
Define a function $u$ on $\overline B(q, \rho_1)$ as follows:
\[
u(x)= \left\{ \begin{array}{ll} \omega(B(q, \rho_1)\setminus\overline G, \partial (B(q,
                       \rho_1)\cap G), x) 
                 & \hbox{ on } x\in B(q, \rho_1)\setminus\overline G,\\
             1 &\hbox{ on } \overline{B(q, \rho_1)\cap G},\\
     0  &\hbox{ on } \partial B(q, \rho_1) \setminus \partial (B(q, \rho_1)\cap G). 
\end{array} \right.
\]
Then $u$ is superharmonic on $B(q, \rho_1)$.
Let $\tilde r_n$ be an increasing sequence going up to $\rho_1$ and $y\in B(q,
\tilde r_1)$. Denote $T_n$ to be the exit time of $X_t$ from $B(q, \tilde
r_n)$. Clearly $T_n$ is increasing and convergent to $T$. Since Brownian
motion is continuous, we have  
\[
\{X_S\in\overline W_{\varepsilon}, S < T\}\subset \underset{n\geq 1}\bigcup 
\{X_{S\wedge T_n}\in\overline W_{\varepsilon}, S < T\},  
\]
where $S\wedge T_n = \min\{S, T_n\}$. So \eqref{4.91} is implied by the 
following:
\begin{equation}\label{4.15}
P^{y}(X_{S\wedge T_n}\in\overline W_{\varepsilon}) = 0, \hbox{ for each } n,
\end{equation}
since by \eqref{4.13}
\begin{eqnarray*}
P^y(X_S\in\overline W_{\varepsilon}) & = &  P^y(X_S\in\overline W_{\varepsilon}, S=T) +
P^y(X_S\in\overline W_{\varepsilon}, S< T)\\ 
& \leq&  P^y(X_T\in\overline W_{\varepsilon}) +P^y(X_S\in\overline W_{\varepsilon}, S<T)\\
&\leq & \lim_{n\to\infty} P^y(X_{S\wedge T_n}\in\overline W_{\varepsilon}, S< T).
\end{eqnarray*}
For a superharmonic function $u$ and for each $n \geq 1$, there exists a 
sequence of increasing superharmonic functions $\{u_j\}$ such that $u_j\in C^2$ and 
\begin{equation}\label{4.16}
\lim_{j\to\infty} u_j = u \hbox{ on }\overline B(q, \tilde r_n)
\end{equation}
(see, for example, Theorem 4.20 of \cite{H}). Applying It\^o's formula to
$u_j(X_{S\wedge T_n})$, we have  
\[
E^y [u_j(X_{S\wedge T_n})] \leq  u_j(y).
\] 
Let $j$ go to infinity and apply \eqref{4.16} to see that
\[
E^y [u(X_{S\wedge T_n})] \leq  u(y).
\]
Consequently
\[
\tau \geq u(y) \geq E^y[u(X_{S\wedge T_n})] \geq P^y(X_{S\wedge T_n}\in
\overline {B(q, \rho_1)\cap G})\geq P^y(X_{S\wedge T_n}\in\overline W_{\varepsilon}).
\]
Letting $\tau \to 0$ proves \eqref{4.15}. 

As a consequence of \eqref{4.9}, since $S$ is upper semi-continuous on 
$\overline B(p, r)\setminus W$, we have almost everywhere with respect to the 
probability measure $P^y$
\begin{equation}\label{4.17}
\limsup_{n\to\infty}S(X_{S_n}) \leq  S(X_{S_\infty})  \hbox{ if }
S_n\uparrow S_{\infty}\leq  T.
\end{equation}   

Let $\eta>0$. Then by \eqref{4.2}, for any $y\in B(q, \rho_1)$, there exists 
$0<\beta = \beta_{y, \eta} < \rho_1-|y-q|$, such that
\begin{equation}\label{4.18}
S(y) - D_{\beta} S(y) \leq a\eta \beta^2.
\end{equation}
Consider a family of stopping times
\[
{\frak S} = \{S\leq  T: S(y) - E^y S(X_S) \leq  a\eta E^y|y - X_S|^2\}.
\]
Define
\[
S_0 = \inf\{ t\geq 0: |X_t- y|\geq \beta_{y,\eta}\}.
\]
Then $X_{S_0}$ is uniformly distributed on $\partial B(y,
\beta_{y,\eta})$. So by \eqref{4.18}, we have
\[
S(y) - E^y S(X_{S_0}) \leq  a\eta E^y|y - X_{S_0}|^2.
\]
Thus $S_0\in {\frak S}$ and hence ${\frak S}$ is not empty.

For a sequence of increasing stopping times $S_n$ in ${\frak S}$, let
$S_\infty = \underset{n\geq 1}\lim S_n$. Then by \eqref{4.17} and
Fatou's lemma,  we have
\[
S(y) - E^y S(X_{S_\infty}) \leq  a\eta E^y |y - X_{S_\infty}|^2.
\]
So $S_\infty\in {\frak S}$. Thus by an argument given in Halmoe [Ha] on page 121\footnote{See [A] for the details.}, there
exists
$S^{*}\in {\frak S}$ such that $S^{*}$ is a maximum of ${\frak S}$. We show that
$S^{*} = T$ almost  everywhere with respect to $P^y$. In fact, if $S^{*} < T$ with positive $P^y$ probability, then    
\[
S^{*}_1 = \inf\{T\geq t\geq S^{*}: |X_t - X_{S^{*}}| \geq \beta_{X_{S^{*}},
\eta}\},
\]
where for $x\in \partial B(q, \rho_1)$, we define $\beta_{x, \eta}$ to be $0$.
Then, clearly, $S^{*}_1 \geq S^{*}$ with strict inequality on 
$\{S^{*} < T\}$. On the other hand, conditional on $X_{S^{*}}$, 
$X_{S^{*}_1}$ is uniformly distributed on the surface of $B(X_{S^{*}}, 
\beta_{X_{S^{*}}, \eta})$, if $S^{*} < S^{*}_1$. So by~\eqref{4.2}
\begin{eqnarray}\label{4.19}
S(X_{S^{*}}) - E^{X_{S^{*}}} S(X_{S^{*}_1}) & =& S(X_{S^{*}}) -
D_{\beta_{X_{S^{*}}, \eta}}S(X_{S^{*}}) \\ \noalign{\vskip4pt} 
& \leq & a\eta |\beta_{X_{S^{*}}, \eta}|^2\nonumber \\ \noalign{\vskip4pt}
& =& a\eta E^{X_{S^*}}|X_{S^*} - X_{S^*_1}|^2. \nonumber
\end{eqnarray}
Hence, by \eqref{4.19}, the strong Markovian property, orthogonality between\break
$X_{S^*_1} - X_{S^*}$ and $X_{S^*} - y$, and $S^*\in {\frak S}$, we have 
\begin{eqnarray*}
S(y)- E^y S(X_{S^*_1}) & =& S(y) - E^y S(X_{S^*}) + E^y[S(X_{S^*}) -
S(X_{S^*_1})]\\ \noalign{\vskip4pt}
& \leq&  a\eta E^y |y - X_{S^*}|^2 + E^y[S(X_{S^*}) - S(X_{S^*_1}), S^* <
S^*_1]\nonumber\\ \noalign{\vskip4pt} 
& =& a\eta E^y |y - X_{S^*}|^2%\nonumber \\ \noalign{\vskip4pt}&&
+\  E^y[S(X_{S^*}) - E^{X_{S^*}}S(X_{S^*_1}),S^* < S^*_1]\nonumber \\ \noalign{\vskip4pt}  
& \leq & a\eta E^y |y - X_{S^*}|^2 + a\eta E^y[E^{X_{S^*}}|X_{S^*} -
X_{S^*_1}|^2, S^* < S^*_1]\nonumber\\ \noalign{\vskip4pt} 
& = &a\eta E^y |y - X_{S^*}|^2 + a\eta E^y |X_{S^*} - X_{S^*_1}|^2\nonumber\\ \noalign{\vskip4pt}
& = & a\eta E^y|y - X_{S^*_1}|^2.\nonumber
\end{eqnarray*}
Thus $S^*_1 \in {\frak S}$. This contradicts the maximality of $S^*$ in
${\frak S}$ since $S^*_1 \geq S^*$ and $S^{*}_1 \ne S^{*}$.  Thus we
have shown that $T\in {\frak S}$. So 
\begin{eqnarray*}
S(y)  - \int_{\partial B(q, \rho_1)} S(z)\, \omega(B(q, \rho_1), dz, y) & = &  S(y) -
E^y S(X_T) \\ \noalign{\vskip4pt}
&\leq&  a\eta E^y |X_T -y|^2 \leq   4a\eta \rho_1^2.
\end{eqnarray*}
This implies \eqref{4.4} by letting $\eta \to 0$. We have finished the proof.
 
\section{Proof of Theorem~\ref{t:4}}\label{s:5}

Without loss of generality, we assume that $q = 0$, the origin. As in the
proof of Theorem~\ref{t:1}, we have that  $\tilde \Delta^{*} f_2(x) =
\tilde\Delta_* f_2(x)$ almost everywhere in $\Bbb  {T}^d
\setminus\{0\}$, and that $\tilde\Delta^* f_2(x)$ is in
$L^1(\Bbb {T}^d\setminus B(0, r))$ for any $r>0$. Consequently,
$\tilde\Delta^*f_2(x) = \tilde\Delta_* f_2(x)$ almost everywhere in
$\Bbb  T^d$. 

As in Section~\ref{s:1}, let $B(x) = \min\{f^*(x), \tilde\Delta^*
f_2(x)\}$. Then by Lemma~\ref{l:1.0}, on 
$\Bbb {T}^d\setminus \{0\}$, $f_*(x) \leq B(x)\leq f^*(x)$ and 
$\tilde\Delta_{*} f_2(x)\leq B(x)\leq \tilde\Delta^{*}f_2(x)$. Thus,
$B(x)\in L^{1}(\Bbb  T^d)$. Consequently, when we proceed  as in 
Section~\ref{s:1}, there exists a function $S^*(x)$, which is harmonic
on $\Bbb {R}^d\setminus M$ and almost everywhere equals
\[
S(x) = f_2(x) + (2\pi)^{-d}\int_{T^d} G(x-y)B(y)\, dy - b_0|x|^2/(2d),
\]
where $M= \{2\mu\pi: \mu\in \Bbb  Z^d\}$. The rest of the proof is
identical to that of Theorem~2 of Shapiro ([Sh, p.479]). 

\medbreak {\it Acknowledgment}.
We are grateful to Jean Bourgain, Robert Kaufman, and Victor Shapiro for helpful conversations.

 \vglue-6pt
\AuthorRefNames [AFR]

\vglue-16pt
\centerline{\ninerm (Received March 17, 1997)}
 \centerline{\ninerm (Revised July 21, 1999)}
 
\end{document}